\documentclass{article}
\usepackage{amsmath,amssymb}
\DeclareOldFontCommand{\sc}{\normalfont\scshape}{\mathrm}
\renewcommand{\Bbb}{\mathbb}
\newtheorem{thm}{Theorem}[section]
\newtheorem{lem}[thm]{Lemma}
\newtheorem{prop}[thm]{Proposition}
\newtheorem{cor}[thm]{Corollary}
\newtheorem{defn}[thm]{Definition}
\newenvironment{defn-new}{\begin{defn} \em}{\end{defn}}
\newtheorem{rem}[thm]{Remark}
\newenvironment{rem-new}{\begin{rem} \em}{\end{rem}}
\newtheorem{ex}[thm]{Example}
\newenvironment{ex-new}{\begin{ex} \em}{\end{ex}}
\newtheorem{prob}[thm]{Problem}
\newenvironment{prob-new}{\begin{prob} \em}{\end{prob}}

\newenvironment{notation-new}{\begin{rem} \em}{\end{rem}}
\newtheorem{agr}[thm]{Agreement}
\newenvironment{agr-new}{\begin{agr} \em}{\end{agr}}

\makeatletter \@addtoreset{equation}{section} \makeatother

\begin{document}

\begin{center}
{\Large {\bf {Trans-Sasakian space forms}}}
\end{center}

\medskip

\renewcommand{\thefootnote}{\fnsymbol{footnote}}

\begin{center}
{\large Prachi$^{1}$ and Mukut Mani Tripathi$^{2}$\footnotemark}
\end{center}

\footnotetext{Corresponding author}

\begin{center}
$^{1}$Department of Mathematics,\\
Institute of Science, \\
Banaras Hindu University,\\
Varanasi 221005, India\\
E-mail: prachisinghsvc203131@gmail.com\\
ORCID: 0009-0006-0218-8109

\medskip

$^{2}$Department of Mathematics,\\
Institute of Science,\\
Banaras Hindu University,\\
Varanasi 221005, India\\
E-mail: mmtripathi66@yahoo.com\\
ORCID: 0000-0002-6113-039X
\end{center}

\begin{quote}
{\bf Abstract.} In this paper, we derive an explicit expression for the Riemann--Christoffel curvature tensor of a trans-Sasakian manifold with pointwise constant $\varphi$-sectional curvature and introduce the notion of trans-Sasakian space forms. As consequences, we obtain explicit formulas for the Ricci tensor and the scalar curvature. We further prove that every $3$-dimensional trans-Sasakian manifold necessarily has pointwise constant $\varphi$-sectional curvature. In addition, we investigate statistical structures on trans-Sasakian manifolds and introduce the notion of trans-Sasakian statistical space forms. We then establish a classification of such statistical structures. Our results show that every trans-Sasakian statistical space form of dimension greater than $3$ is necessarily almost-trivial. In the $3$-dimensional case, we establish a rigidity theorem for trans-Sasakian statistical space forms and provide explicit examples showing the existence of non-trivial $3$-dimensional trans-Sasakian statistical space forms.
\end{quote}


\medskip

\noindent {\em $2020$ Mathematics Subject Classification\/}{\rm :} {53B12, 53C25, 53D15}


\section{Introduction}
Curvature plays a central role in the study of the geometric properties of manifolds. In Riemannian geometry, the sectional curvatures determine the Riemann--Christoffel curvature tensor of a Riemannian manifold, while in Kaehler geometry the holomorphic sectional curvatures completely characterize the Riemann--Christoffel curvature tensor. In contact metric geometry, Moskal showed that the Riemann--Christoffel curvature tensor of a Sasakian manifold is characterized by its $\varphi$-sectional curvatures completely \cite[Moskal 1966]{Moskal-66-Thesis}. Subsequently, Ogiue derived an explicit expression for the curvature tensor of Sasakian manifolds with (pointwise) constant $\varphi$-sectional curvature \cite[Ogiue 1964]{Ogiue-64-Kodai}. Such manifolds are called Sasakian space forms \cite[Harada 1972]{Harada-72-TGU}. Analogous curvature formulas have since been obtained for several classes of almost contact metric manifolds having (pointwise) constant $\varphi$-sectional curvature, including Kenmotsu space forms \cite[Kenmotsu 1972]{Kenmotsu-72-Tohoku} and cosymplectic space forms \cite[Ludden 1970]{Ludden-70-JDG}. 

In 1985, Oubina introduced trans-Sasakian manifolds as a special class of almost contact metric manifolds \cite[Oubina 1985]{Oubina-85-PMD}. Subsequently, in 1990, Blair and Oubina established a necessary and sufficient condition for an almost contact metric manifold to be a trans-Sasakian manifold of type $(\alpha,\beta)$ \cite[Blair and Oubina 1990]{Blair-Oubina-90-PM}. This class contains several important subclasses, including cosymplectic manifolds \cite[Blair 1967]{Blair-67-JDG}, $\alpha$-Sasakian manifolds \cite[Janssens and Vanhecke 1981]{Janssens-Vanhecke-81-KMJ}, Sasakian manifolds \cite[Sasaki and Hatakeyama 1962]{Sasaki-Hatakeyama-62-JMSJ}, $\beta$-Kenmotsu manifolds \cite[Janssens and Vanhecke 1981]{Janssens-Vanhecke-81-KMJ}, Kenmotsu manifolds \cite[Kenmotsu 1972]{Kenmotsu-72-Tohoku}. In 1992, Marrero characterized the local structure of trans-Sasakian manifolds completely \cite[Marrero 1992]{Marrero-92-AMPA}. Although many geometric properties of trans-Sasakian manifolds have been extensively studied, an explicit expression for the Riemann--Christoffel curvature tensor under the assumption of pointwise constant $\varphi$-sectional curvature has not been previously available. 

The primary goal of this paper is to derive such a formula, thereby resolving a longstanding problem in the theory of trans-Sasakian manifolds that has remained open since their introduction by Oubina in 1985. Motivated by this characterization, we introduce the notion of a trans-Sasakian space form. This provides a complete curvature characterization in the trans-Sasakian setting and generalizes several classical curvature formulas established for its important subclasses, including Sasakian, Kenmotsu, cosymplectic, $\alpha$-Sasakian, and $\beta$-Kenmotsu manifolds.

This result provides a powerful tool for further investigations, including classification results, curvature characterizations, inequalities, curvature of product trans-Sasakian space forms, the study of geometric and analytic properties of trans-Sasakian manifolds etc. Consequently, it contributes significantly to a deeper and more comprehensive understanding of the geometry of trans-Sasakian manifolds and open several new directions for further research.

Motivated by the above curvature characterization, it is natural to investigate its analogue in the statistical setting. The concept of a statistical manifold was first introduced by Amari \cite[Amari 1985]{Amari-85-Book} as a differential geometric framework for the study of families of probability distributions. Since then, statistical geometry has developed strong connection with differential geometry (see \cite[Matsuzoe 2010]{Matsuzoe-10-Chap}) as well as Hessian geometry (see \cite[Shima 2013]{Shima-13-Chapter}). Owing to its ability to describe statistical models geometrically, the theory of statistical manifolds has become a fundamental tool in information geometry, with important applications in machine learning, neural networks, statistical inference, optimization, and mathematical physics (see \cite[Amari and Kawanabe 1994]{Amari-Kawanabe-94-SK}, \cite[Sun and Marchand 2014]{Sun-Marchand-14-Conf}). 

In addition to its theoretical significance, statistical geometry has found numerous applications in information geometry. For example, in \cite[Ali {\em et al.} 2010]{Ali-Cafaro-Kim-Mancini-10-PhyA}, the authors analytically studied the asymptotic temporal behavior of the information geometric complexity (IGC) of finite-dimensional Gaussian statistical manifolds in the presence of microcorrelations among microvariables. Subsequently, \cite[Gomez 2017]{Gomez-17-PhyA} employed information-geometric techniques to extend the classical ergodic hierarchy by introducing geometric analogues of the ergodic, mixing, and Bernoulli levels for statistical models defined on curved manifolds. 

These developments demonstrate that curvature is one of the fundamental geometric invariants in information geometry, as it reflects the intrinsic geometry and complexity of statistical models. Consequently, introducing new classes of statistical manifolds together with explicit curvature characterizations enlarges the geometric framework available for information geometry and provides a foundation for future investigations of statistical modeling, information-geometric complexity, and related problems in mathematical physics.

Meanwhile, almost contact manifolds form an important class of odd-dimensional manifolds in Riemannian geometry and have been extensively studied due to their rich geometric structure and their applications in the study of dynamical systems and mathematical physics. Combining statistical and almost contact structures leads to a natural geometric framework that unifies probabilistic and differential geometric methods, thereby providing new tools for investigating geometric models arising in both statistics and physics.

It is worth noting that every Riemannian structure can be regarded as a particular statistical structure for which the difference tensor field ${\sc B}^{\sc d}$ vanishes. Consequently, many geometric structures defined in the Riemannian setting admit natural generalizations to the statistical framework. Motivated by this observation, Furuhata {\em et al.} \cite[Furuhata {\em et al.} 2017]{Furuhata-Hasegawa-Okuyama-Sato-Shahid-17-JGP} introduced the notion of a Sasakian statistical manifold and established conditions under which a real hypersurface of a holomorphic statistical manifold carries such a structure. In another work, Furuhata {\em et al.} \cite[Furu-Hasegawa-Okuyama-Sato-17-JG]{Furu-Hasegawa-Okuyama-Sato-17-JG} developed the theory of Kenmotsu statistical manifolds by investigating statistical structures compatible with Kenmotsu geometry. Subsequently, Kazan \cite[Kazan 2019]{Kazan-19-PhyA} initiated the study of trans-Sasakian statistical manifolds by establishing several fundamental curvature identities and investigating conformally-projectively flat trans-Sasakian statistical manifolds. However, the theory of trans-Sasakian statistical space forms and the corresponding statistical curvature characterization have remained unexplored.

Motivated by the above developments, the second major objective of this paper is to investigate the statistical structures on trans-Sasakian space forms and to introduce the notion of trans-Sasakian statistical space forms.

The structure of the paper is organized as follows. Section~\ref{Sect-Prel} reviews the fundamental definitions and preliminary results on statistical manifolds and trans-Sasakian manifolds. In Section~\ref{Sect-Trans-SasSF}, we establish a necessary and sufficient condition, expressed in terms of the Riemann--Christoffel curvature tensor, for a trans-Sasakian manifold to admit pointwise constant $\varphi$-sectional curvature. Corresponding formulas for the Ricci tensor and the scalar curvature are also derived. Furthermore, we show that every $3$-dimensional trans-Sasakian manifold has pointwise constant $\varphi$-sectional curvature, and we construct an explicit example verifying the main results. In Section \ref{sec-trans-sasa-stat}, we establish several results concerning trans-Sasakian statistical manifolds. Section \ref{sec-tran-Sasa-Stat-nonproper} shows that every trans-Sasakian statistical space form of dimension greater than $3$ is necessarily almost-trivial, and an explicit example is provided. In Section \ref{sec-3D-trans-sasa-stat}, we establish a rigidity theorem for $3$-dimensional trans-Sasakian statistical space forms and provide explicit examples showing the existence of non-trivial $3$-dimensional trans-Sasakian statistical space forms.

\section{Preliminaries\label{Sect-Prel}}

\subsection{Statistical manifolds}

Let $(M,g)$ be a Riemannian manifold and $\nabla$ denote the Levi--Civita connection of $g$. Suppose that ${\sc D}$ is an affine connection on $M$, then the dual connection ${\sc D}^{\ast}$ is defined by 
\begin{equation}
Xg\left(Y,Z\right) =g\left({\sc D}_{X}Y,Z\right) +g\left(Y,{\sc D}_{X}^{\ast}Z\right)  \label{eq-nabla-star}
\end{equation}
for all vector fields $X$, $Y$, $Z$ on $M$ (see \cite[Amari 1985]{Amari-85-Book}).

A pair $({\sc D},g)$ is said to be a statistical structure or a Codazzi structure if ${\sc D}$ is torsion-free and
\begin{equation}
\left({\sc D}_{X}g\right) \left(Y,Z\right) = \left({\sc D}_{Y}g\right) \left(X,Z\right)  \label{eq-Statis-struc}
\end{equation}
for all vector fields $X$, $Y$, $Z$ on $M$ (\cite[Furuhata 2002]{Furuhata-2002-Chap}, \cite[Matsuzoe 2010]{Matsuzoe-10-Chap}). Then $(M,{\sc D},g)$ is called a statistical manifold or a Codazzi manifold. Furthermore, when both ${\sc D}$ and ${\sc D}^{\ast}$ are torsion-free, one has \cite[Matsuzoe 2010]{Matsuzoe-10-Chap} 
\begin{equation}
{\sc D} +{\sc D}^{\ast}=2\nabla.  \label{eq-nabla-star-zero}
\end{equation}
If $({\sc D},g)$ is a statistical structure on $M$, then its dual pair $({\sc D}^{\ast},g)$ also defines a statistical structure on $M$ \cite[Matsuzoe 2010]{Matsuzoe-10-Chap}. A Riemannian structure $(\nabla,g)$ naturally forms a statistical structure. Such a structure is referred to as a Riemannian statistical structure or a trivial statistical structure \cite[Furuhata {\em et al.} 2017]{Furuhata-Hasegawa-Okuyama-Sato-Shahid-17-JGP}.

The difference tensor field ${\sc B}^{\sc d}$ on $(M,g,{\sc D})$ is given by \cite[Furuhata 2009]{Furuhata-09-DGA}
\begin{equation}
{\sc B}^{\sc d}_{X}Y = {\sc D}_{X}Y - \nabla_{X}Y,  \label{eq-K(X,Y)}
\end{equation}
which satisfies
\begin{equation}
{\sc B}^{\sc d}_{X}Y = {\sc B}^{\sc d}_{Y}X,\quad g\left({\sc B}^{\sc d}_{X}Y,Z\right) = g\left({\sc B}^{\sc d}_{X}Z,Y\right).
\label{eq-K-properties}
\end{equation}
for all vector fields $X$, $Y$, $Z$ on $M$. If the tensor field ${\sc B}^{\sc d}$ satisfies (\ref{eq-K-properties}), then $(M,\nabla+{\sc B}^{\sc d},g)$ defines a statistical manifold \cite[Furuhata {\em et al.} 2017]{Furuhata-Hasegawa-Okuyama-Sato-Shahid-17-JGP}. If the difference tensor field ${\sc B}^{\sc d}$ vanishes then the statistical structure is called a trivial statistical structure.

\begin{defn-new}
\cite[Furuhata and Hasegawa 2016]{Furuhata-Hasegawa-16-Chap} Let $(M,g,{\sc D})$ be a statistical manifold and let ${\sc D}^{\ast}$ be the dual connection of ${\sc D}$. Let $R$ and $R^{\ast}$ be the curvature tensor corresponding to ${\sc D}$ and ${\sc D}^{\ast}$, define by
\begin{eqnarray}
R^{\sc d}(X,Y)Z &=&{\sc D}_{X}{\sc D}_{Y}Z-{\sc D}_{Y}{\sc D}_{X}Z - {\sc D}_{[X,Y]}Z,  \label{eq-R-nabla} \\
R^{{\sc d}^{\ast}}(X,Y)Z &=& {\sc D}_{X}^{\ast}{\sc D}_{Y}^{\ast}Z-{\sc D}_{Y}^{\ast}{\sc D}_{X}^{\ast}Z-{\sc D}_{[X,Y]}^{\ast}Z \label{eq-R-star-nabla}
\end{eqnarray}
for all vector fields $X$, $Y$, $Z$ on $M$.
Define
\begin{equation}
R^{\sc s} = \frac{1}{2}(R^{\sc d} + R^{{\sc d}^{\ast}}),  \label{eq-S-statis-cur}
\end{equation}
then $R^{\sc s}$ is called the statistical curvature tensor.
\end{defn-new}

Let $(M,{\sc D},g)$ be a statistical manifold, and $\Pi_{p}^{2}$ be the $2$-plane section in $T_{p}M$ spanned by the vector fields $X$ and $Y$. Then \begin{equation}
K^{\sc s}(X,Y)=\frac{g\left(R^{\sc s}\left(X,Y\right) Y,X\right)}{g\left(X,X\right)g\left(Y,Y\right) - g\left(X,Y\right)^{2}}  \label{eq-statis-sec-cur}
\end{equation}
is called the statistical sectional curvature (see \cite[Furuhata and Hasegawa 2016]{Furuhata-Hasegawa-16-Chap}).

If $R$ is the curvature tensor associated with the Levi--Civita connection $\nabla$. Then (see {\rm \cite[Furuhata {\em et al.} 2017]{Furuhata-Hasegawa-Okuyama-Sato-Shahid-17-JGP}})
\begin{equation}
R^{\sc s}\left(X,Y\right) Z = R\left(X,Y\right) Z+[{\sc B}^{\sc d}_{X},{\sc B}^{\sc d}_{Y}]Z
\label{eq-statistical-cur}
\end{equation}
for all vector fields $X$, $Y$, $Z$ on $M$, where
\begin{equation}
[{\sc B}^{\sc d}_{X},{\sc B}^{\sc d}_{Y}] = {\sc B}^{\sc d}_{X}{\sc B}^{\sc d}_{Y}-{\sc B}^{\sc d}_{Y}{\sc B}^{\sc d}_{X}.  \label{eq-commu-K}
\end{equation}

\subsection{Trans-Sasakian manifolds}
Let $M$ be a smooth manifold of dimension $(2n+1)$. An almost contact metric structure on $M$ determined by a quadruple $(\varphi,\xi,\eta,g)$, where $\varphi$ is a $(1,1)$-type tensor field, $\xi$ is a vector field, $\eta$ is a $1$-form, and $g$ is a compatible Riemannian metric satisfying
\[
{\varphi}^{2} = -I+\eta \otimes \xi , \quad \eta(\xi) = 1,\quad \varphi \xi = 0, \quad \eta \circ \varphi = 0, 
\]
together with
\begin{equation}
g\left(\varphi {X},\varphi {Y}\right) = g\left(X,Y\right) - \eta \left(X\right) \eta \left(Y\right) \label{eq-g(phiX,phiY)}
\end{equation}
for all vector fields $X$, $Y$ on $M$. A manifold $M$ equipped with the structure satisfying above conditions is called an almost contact metric manifold. For details, we refer to (\cite[Sasaki 1960]{Sasaki-60-TMJ}, \cite[Blair 2010]{Blair-10-Book}).

\begin{defn-new}
\cite[Oubina 1985]{Oubina-85-PMD} An almost contact metric manifold $(M,\varphi,\xi,\eta,g)$ is called trans-Sasakian if the associated almost Hermitian manifold $(M\times {\Bbb R},J^{\times},g^{\times})$ belongs to the Gray--Hervella  class ${\cal W}_{4}$ \cite[Gray and Hervella 1980]{Gray-Hervella-80-AMPA}, here the almost complex structure $J^{\times}$ is defined by 
\[
J^{\times}\left(X,f\frac{d}{dt}\right) = \left(\varphi{X}-f\xi,\eta (X)\frac{d}{dt}\right) 
\]
for any vector field $X$ on $M$ and any smooth function $f$ on $M\times {\Bbb R}$ while $g^{\times}$ represents the product metric on $M\times {\Bbb R}$.
\end{defn-new}

It was shown in \cite[Blair and Oubina 1990]{Blair-Oubina-90-PM}, an almost contact metric manifold $(M,\varphi,\xi,\eta,g)$ admits a trans-Sasakian structure if and only if
\begin{equation}
\left(\nabla_{X}\varphi \right) Y=\alpha \left(g\left(X,Y\right) \xi - \eta \left(Y\right) X\right) + 
\beta \left(g\left(\varphi{X},{Y}\right)\xi - \eta \left(Y\right) \varphi {X}\right)  
\label{eq-tr-sas-def-a}
\end{equation}
for some smooth functions $\alpha$ and $\beta$ on $M$, where $\nabla$ is the Levi-Civita connection of $g$. Such a manifold is then said to be a trans-Sasakian manifold of type $(\alpha,\beta)$ \cite[Blair and Oubina 1990]{Blair-Oubina-90-PM}. 

\begin{rem-new}\label{rem-subcls-tsm}
An $(\alpha,\beta)$ trans-Sasakian structure reduces to
\begin{itemize}
\item[{\bf (a)}] ${\cal C}_{6}$-structure (\cite[Chinea and Gonzalez 1990]{Chinea-Gonzalez-90-AMPA}, \cite[Marrero 1992]{Marrero-92-AMPA}) if $\beta=0$,

\item[{\bf (b)}] an $\alpha$-Sasakian structure \cite[Janssens and Vanhecke 1981, Theorem 2.6]{Janssens-Vanhecke-81-KMJ} if $\beta=0\neq \alpha =$ constant,
    
\item[{\bf (c)}] a Sasakian structure \cite[Blair 2010]{Blair-10-Book} if $\beta=0$, $\alpha=1$,

\item[{\bf (d)}] ${\cal C}_{5}$-structure (\cite[Chinea and Gonzalez 1990]{Chinea-Gonzalez-90-AMPA}, \cite[Marrero 1992]{Marrero-92-AMPA}) if $\alpha=0$,

\item[{\bf (e)}] a $\beta$-Kenmotsu structure \cite[Janssens and Vanhecke 1981, Theorem 2.8]{Janssens-Vanhecke-81-KMJ} if $\alpha =0\neq \beta =$ constant,

\item[{\bf (f)}] a Kenmotsu structure \cite[Kenmotsu 1972]{Kenmotsu-72-Tohoku} if $\alpha=0$, $\beta=1$,

\item[{\bf (g)}] a cosymplectic structure \cite[Blair 1967]{Blair-67-JDG} if $\alpha =0=\beta$.
\end{itemize}
\end{rem-new}

For a trans-Sasakian manifold $(M,\varphi,\xi,\eta,g)$, the following identities hold \cite[De and Tripathi 2003]{DeUC-Tri-03-KMJ}:
\begin{eqnarray}
R\left(X,Y\right) \xi &=&\left({\alpha}^{2} - {\beta}^{2}\right) \left(\eta\left(Y\right)X - \eta\left(X\right)Y\right) + 2\alpha\beta\left(\eta \left(Y\right)\varphi{X}-\eta\left(X\right) \varphi{Y}\right) \nonumber \\
&&+ {\rm d}\alpha \left(Y\right) \varphi {X} - {\rm d}\alpha \left(X\right)\varphi {Y} + {\rm d}\beta \left(Y\right){\varphi}^{2}X - {\rm d}\beta \left(X\right) {\varphi}^{2}Y,\label{eq-tr-sasa-cur-ope}
\end{eqnarray}
\begin{equation}
2\alpha\beta + {\rm d}\alpha\left(\xi\right) = 0, \label{eq-condition-alpha-beta-xi}
\end{equation}
\begin{equation}
\nabla_{X}\xi = -\alpha \varphi{X}+\beta \left(X-\eta \left(X\right) \xi \right), \label{eq-nabla_X(xi)}
\end{equation}
\begin{equation}
\left(\nabla_{X}\eta \right) Y = \alpha g\left(X,\varphi {Y}\right) +\beta \left(g\left(X,Y\right) - \eta \left(X\right) \eta \left(Y\right) \right) \label{eq-nabla_X(eta)}
\end{equation}
for all vector fields $X$, $Y$ on $M$, where $R$ is the Riemann curvature tensor.

\begin{defn-new}
The symmetric product $\Theta_{1}\circledcirc \Theta_{2}$ and the Kulkarni--Nomizu product 
$\Theta_{1}\circledast \Theta_{2}$ of any $(0,2)$-tensors $\Theta_{1}$ and $\Theta_{2}$ on a 
smooth manifold $M$ are $(0,4)$-tensors defined as follows (see \cite[Tripathi 2025, Eqs.~(4.1), 
(3.1)]{Tri-25-CM})
\begin{equation}
\Theta_{1}\circledcirc \Theta_{2}=\Theta_{1}\otimes \Theta_{2}+\Theta_{2}\otimes \Theta_{1},  \label{eq-V-Sym-Product}
\end{equation}
\begin{eqnarray}
(\Theta_{1}\circledast \Theta_{2})(X,Y,Z,W) &=&\Theta_{1}(Y,Z)\Theta_{2}(X,W)-\Theta_{1}(X,Z)\Theta_{2}(Y,W)  \nonumber \\
&&+\Theta_{2}(Y,Z)\Theta_{1}(X,W)-\Theta_{2}(X,Z)\Theta_{1}(Y,W)
\label{eq-V-KN-Product}
\end{eqnarray}
for all vector fields $X$, $Y$, $Z$, $W$ on $M$.
\end{defn-new}

A $(0,4)$-tensor $T$ on a smooth manifold $M$ is said to be a Kulkarni--Nomizu tensor if (see 
\cite[Kulkarni 1969]{Kulkarni-69-BAMS}, \cite[Tripathi 2025, Lemma 5.1]{Tri-25-CM}) 
\begin{equation}
T(X,Y,Z,W) = -T(Y,X,Z,W),  \label{KN-tensor-1}
\end{equation}
\begin{equation}
T(X,Y,Z,W) = T(Z,W,X,Y),  \label{KN-tensor-3}
\end{equation}
\begin{equation}
T(X,Y,Z,W)+T(Y,Z,X,W)+T(Z,X,Y,W) = 0  \label{KN-tensor-4}
\end{equation}
for all vector fields $X$, $Y$, $Z$, $W$ on $M$.

\section{Trans-Sasakian space forms\label{Sect-Trans-SasSF}}
The aim of this section is to derive an explicit expression for the Riemann--Christoffel curvature tensor of trans-Sasakian manifolds with pointwise constant $\varphi$-sectional curvature. To this end, we first establish several preliminary results that will be used in the proof of the main theorem.

\begin{lem}
\label{lem-sstf-kntf} Let $\Theta_{1}$ and $\Theta_{2}$ be two skew symmetric tensors of type $(0,2)$ on a smooth manifold $M$. Then, the $\left(0,4\right)$-tensor $T$ defined by
\begin{equation}
T = \Theta_{1}\circledast \Theta_{2} - 2\Theta_{1}\circledcirc \Theta_{2}, \label{eq-Theta1-Theta2-knp-snp}
\end{equation}
is a Kulkarni--Nomizu tensor on $M$. Consequently, if $\Theta_{1}\circledast \Theta_{2}$ is a 
Kulkarni--Nomizu tensor, then $\Theta_{1}\circledcirc \Theta_{2}$ is also a Kulkarni--Nomizu 
tensor. Moreover,
\begin{equation}
\Theta_{1}\circledcirc \Theta_{2} = -\Theta_{1}\circledast \Theta_{2}. \label{eq-Theta1-smp-Theta2}
\end{equation}
\end{lem}

\noindent {\bf Proof.} It is straightforward to verify that the tensor $T$, defined by 
(\ref{eq-Theta1-Theta2-knp-snp}), satisfies conditions (\ref{KN-tensor-1}), (\ref{KN-tensor-3}) 
and (\ref{KN-tensor-4}), hence $T$ is a Kulkarni--Nomizu tensor. If $\Theta_{1}\circledast 
\Theta_{2}$ is a Kulkarni--Nomizu tensor then from (\ref{eq-Theta1-Theta2-knp-snp}), 
$\Theta_{1}\circledcirc \Theta_{2}$ is also a Kulkarni--Nomizu tensor. Next, applying the Bianchi 
identity (\ref{KN-tensor-4}) to Kulkarni--Nomizu tensor $\Theta_{1}\circledcirc \Theta_{2}$, we 
obtain (\ref{eq-Theta1-smp-Theta2}). $\blacksquare$

On a Riemannian manifold $(M,g)$, the scalar curvature $\tau$ is given by
\begin{equation}
\tau=\frac{1}{2}{\rm tr}(S), \label{eq-scal-cur}
\end{equation}
where $S$ is the Ricci tensor of $M$.

According to \cite[Marrero 1992]{Marrero-92-AMPA}, a trans-Sasakian manifold of dimension $\geq 5$ is necessarily either cosymplectic, an $\alpha$-Sasakian or a $\beta$-Kenmotsu manifold. As a consequence, proper trans-Sasakian manifolds occur only in the 3-dimensional case. In a trans-Sasakian manifold of dimension $3$, the Riemann curvature tensor $R$ and the scalar curvature $\tau$ are related by (cf.~\cite[De and Tripathi 2003, Eq.~(20)]{DeUC-Tri-03-KMJ}) 
\begin{eqnarray}
R &=&\left(2\tau +{\rm d}\beta \left(\xi\right) - \left({\alpha}^{2} - {\beta}^{2}\right) \right) g\circledast g - \left(\tau + {\rm d}\beta \left(\xi\right) - 3\left({\alpha}^{2} - {\beta}^{2}\right) \right) g\circledast\left(\eta \otimes \eta \right)  \nonumber\\
&& + g\circledast \left(\eta \otimes \left(\varphi{\rm grad}\left(\alpha\right)\right)^{\sharp} + \left(\varphi{\rm grad}\left(\alpha\right)\right)^{\sharp}\otimes \eta \right) - g\circledast\left(\eta \otimes{\rm d}\beta +{\rm d}\beta \otimes \eta \right),\label{eq-3D-tsm-R(X,Y,Z)} 
\end{eqnarray}
where 
\begin{equation}
\left(\varphi {\rm grad}\left(\alpha \right) \right)^{\sharp}\left(\cdot \right) = g\left(\cdot,\varphi {\rm grad}\left(\alpha \right) \right) .  \label{eq-grad(sharp)}
\end{equation}

\begin{prop}
In a smooth manifold $M$, let $(\varphi,\xi,\eta,g)$ be a trans-Sasakian structure of type $(\alpha,\beta)$. Then,
\begin{equation}
{\rm d}\alpha(\varphi{Y}) {\varphi}^{2}X-{\rm d}\alpha(\varphi{X}) {\varphi}^{2}Y = ({\rm d}\alpha(Y) -{\rm d}\alpha(\xi) \eta(Y)) \varphi {X}- ({\rm d}\alpha(X) -{\rm d}\alpha(\xi)\eta(X)) \varphi {Y}  \label{eq-3D-tsm-equi-R(X,Y)xi}
\end{equation}
for all vector fields $X$, $Y$ on $M$.
\end{prop}

\noindent {\bf Proof.} If ${\rm dim}(M)=3$, taking (\ref{eq-3D-tsm-R(X,Y,Z)}) into account, it follows that
\begin{eqnarray}
R\left(X,Y\right) \xi &=&\left({\alpha}^{2}-{\beta}^{2}\right) \left(\eta \left(Y\right) X-\eta \left(X\right) Y\right) + {\rm d}\beta \left(Y\right) {\varphi}^{2}X-{\rm d}\beta \left(X\right) {\varphi}^{2}Y  \nonumber\\
&&-{\rm d}\alpha \left(\varphi{X}\right) {\varphi}^{2}Y+{\rm d}\alpha \left(\varphi{Y}\right) 
{\varphi}^{2}X  \label{eq-3D-tsm-R(Y,Y)xi}
\end{eqnarray}
for all vector fields $X$, $Y$ on $M$. Comparing (\ref{eq-tr-sasa-cur-ope}) and (\ref{eq-3D-tsm-R(Y,Y)xi}), and using (\ref{eq-condition-alpha-beta-xi}), we obtain (\ref{eq-3D-tsm-equi-R(X,Y)xi}). If ${\rm dim}(M)>3$, then $\alpha$ reduces to be a constant and (\ref{eq-3D-tsm-equi-R(X,Y)xi}) holds identically. $\blacksquare$

The well-known Sasaki $2$-form ${\varphi}^{\flat}$, associated with an almost contact metric structure $(\varphi,\xi,\eta,g)$ in a smooth manifold $M$, is given by 
\begin{equation}
{\varphi}^{\flat}\left(X,Y\right) = g\left(X,\varphi Y\right) \label{eq-funda-2-form}
\end{equation}
for all vector fields $X$, $Y$ on $M$.

\begin{prop}
In a smooth manifold $M$, let $\left(\varphi,\xi,\eta,g\right)$ be a trans-Sasakian structure of type $(\alpha,\beta)$. Then, 
\begin{eqnarray}
&&R\left(X,Y,Z,\varphi{W}\right)+R\left(X,Y,\varphi{Z},W\right)  \nonumber\\
&&\;=\left\{ g\circledast \left({\rm d}\alpha \otimes \eta \right) -{\varphi}^{\flat}\circledast \left({\rm d}\beta \otimes \eta \right) + \left(\alpha^{2}-{\beta}^{2}\right) g\circledast {\varphi}^{\flat} +\alpha \beta \left(g\circledast g-{\varphi}^{\flat} \circledast {\varphi}^{\flat} \right) \right\} \left(X,Y,Z,W\right) \label{eq-R(X,Y)PhiZ-Phi(R(X,Y)Z)}
\end{eqnarray}
for all vector fields $X$, $Y$, $Z$, $W$ on $M$.
\end{prop}

\noindent {\bf Proof.} Ricci identity for the tensor $\varphi$ of type $(1,1)$ yields 
\begin{eqnarray}
g(R(X,Y) \varphi Z,W) - g(\varphi R(X,Y) Z,W) &=& g(\nabla_{X}((\nabla_{Y}\varphi) Z),W) - g((\nabla_{Y}\varphi)(\nabla_{X}Z),W) \nonumber \\
&& - g(\nabla_{Y}((\nabla_{X}\varphi) Z),W) + g((\nabla_{X}\varphi) (\nabla_{Y}Z),W) \nonumber \\
&& - g((\nabla_{[X,Y]}\varphi)Z,W)\label{eq-trsm-Ricci-identity-0}
\end{eqnarray}
for all vector fields $X$, $Y$, $Z$, $W$ on $M$.
By making use of (\ref{eq-tr-sas-def-a}), (\ref{eq-nabla_X(xi)}), (\ref{eq-nabla_X(eta)}), (\ref{eq-V-KN-Product}) and (\ref{eq-funda-2-form}) in (\ref{eq-trsm-Ricci-identity-0}), we arrive at 
\begin{eqnarray}
&&g(R(X,Y)\varphi Z,W) - g(\varphi R(X,Y)Z,W)   \nonumber \\
&&\;\;=\left(g\circledast({\rm d}\alpha \otimes \eta) - \varphi^{\flat}\circledast ({\rm d}\beta \otimes \eta) + ({\alpha}^{2} - {\beta}^{2}) g\circledast \varphi^{\flat} \right. \nonumber\\
&&\quad \qquad \left. + \alpha\beta(g\circledast g-\varphi^{\flat}\circledast \varphi^{\flat})\right)({X,Y,Z,W}).\label{eq-tps-Riem-cur-13a}
\end{eqnarray}
By using $g(\varphi R(X,Y)Z,W)=-g(R(X,Y)Z,\varphi W)$ in (\ref{eq-tps-Riem-cur-13a}), we obtain (\ref{eq-R(X,Y)PhiZ-Phi(R(X,Y)Z)}).
$\blacksquare$

Consider a $(2n+1)$-dimensional almost contact metric manifold $(M,\varphi,\xi,\eta,g)$. A $2$-plane section in $T_{p}M$ is known to be a $\varphi$-section if there is a vector $X\in T_{p}M$ orthogonal to $\xi$ such that $\{X,\varphi {X}\}$ spans the section. The sectional curvature $K(X,\varphi {X})$, denoted by $c(X)$, is called a $\varphi$-sectional curvature, given by (cf. \cite[Ogiue 1964]{Ogiue-64-Kodai}) 
\begin{equation}
c(X) = K(X,\varphi {X}) = \frac{R(X,\varphi {X},\varphi{X},X)}{g(X,X)g(\varphi{X},\varphi {X}) - g(X,\varphi {X})^{2}} = \frac{R(X,\varphi {X},\varphi {X},X)}{g(X,X)^{2}}.  \label{eq-phi-para-holo-sec-cur}
\end{equation}
We now present the central result of this study.
\begin{thm}
\label{th-ts-sf} In a smooth manifold $M$, let $(\varphi,\xi,\eta,g)$ be a trans-Sasakian structure of type $(\alpha,\beta )$. For $M$, the necessary and sufficient condition of having pointwise constant $\varphi$-sectional curvature $c$ is that its Riemannian-Christoffel curvature tensor is given by 
\begin{eqnarray}
R &=&\frac{c+3\left({\alpha}^{2}-{\beta}^{2}\right)}{8}g\circledast g+\frac{c-\left({\alpha}^{2}-{\beta}^{2}\right)}{4}\left(\frac{1}{2}{\varphi}^{\flat} \circledast {\varphi}^{\flat} -{\varphi}^{\flat} \circledcirc {\varphi}^{\flat} \right)  \nonumber \\
&&+\left({\rm d}\beta \left(\xi \right) -\frac{c-\left({\alpha}^{2}-\beta ^{2}\right)}{4}\right) g\circledast \left(\eta \otimes \eta \right)\nonumber\\
&&+g\circledast \left(\left(\varphi {\rm grad}\left(\alpha \right)\right)^{\sharp}\otimes \eta +\eta \otimes \left(\varphi {\rm grad}\left(\alpha \right) \right)^{\sharp}-\left({\rm d}\beta \otimes \eta +\eta \otimes d\beta \right) \right).  \label{eq-tps-Riem-cur-for-33}
\end{eqnarray}
\end{thm}

\noindent {\bf Proof.} From (\ref{eq-phi-para-holo-sec-cur}), we obtain
\begin{equation}
R(X,\varphi{X},\varphi{X},X) = cg(X,X)^{2}.  
\label{eq-phi-sec-curv}
\end{equation}
Replacing $W$ by $\varphi {W}$ in (\ref{eq-R(X,Y)PhiZ-Phi(R(X,Y)Z)}), then interchanging $X$ with $Z$, and $Y$ with $W$ in the resulting equation, we derive
\begin{eqnarray}
&&R\left(\varphi{X},\varphi{Y},Z,W\right) +R\left(X,{\varphi}^{2}Y,Z,W\right) \nonumber \\
&&\quad =\left\{ g\circledast \left({\rm d}\alpha \otimes \eta \right)-{\varphi}^{\flat} \circledast \left({\rm d}\beta \otimes \eta \right) +\left(\alpha^{2}-{\beta}^{2}\right) g\circledast {\varphi}^{\flat} \right.  \nonumber \\
&&\qquad \quad \left. +\alpha \beta \left(g\circledast g-{\varphi}^{\flat} \circledast {\varphi}^{\flat} \right) \right\} 
\left(Z,W,X,\varphi{Y}\right). \label{eq-R(phiX,phiY)Z+R(X,Y)Z-0}
\end{eqnarray}
Let ${\cal D}=\{X\in TM:\eta \left(X\right)=0\}$ denote the contact distribution on $M$. Then every vector field $X$ on $M$ admits a unique decomposition of the form
\begin{equation}
X=X_{\cal D}+\eta \left(X\right) \xi, \quad X_{\cal D}\in {\cal D}.  \label{eq-barX}
\end{equation}
Replacing $X$ and $Z$ by $X_{\cal D}$; $Y$ by $\varphi X_{\cal D}$; and $W$ by 
$Y_{\cal D}$ in (\ref{eq-R(X,Y)PhiZ-Phi(R(X,Y)Z)}), we get
\begin{equation}
R\left(X_{\cal D},\varphi X_{\cal D},\varphi X_{\cal D},Y_{\cal D}\right) =R\left(X_{\cal D},\varphi 
X_{\cal D},\varphi Y_{\cal D},X_{\cal D}\right). \label{eq-tps-Riem-cur-for-2}
\end{equation}
Replacing $Y$ by $\varphi Y_{\cal D}$; $Z$ by $Y_{\cal D}$; $W$ and $X$ by $X_{\cal D}$ in 
(\ref{eq-R(X,Y)PhiZ-Phi(R(X,Y)Z)}), we get
\begin{eqnarray}
&&R\left(X_{\cal D},\varphi Y_{\cal D},\varphi Y_{\cal D}, X_{\cal D}\right) +R\left(X_{\cal D},\varphi Y_{\cal D},Y_{\cal D},\varphi X_{\cal D}\right)   \nonumber \\
&&\qquad =\left({\alpha}^{2}-{\beta}^{2}\right) \left\{-g\left(X_{\cal D},Y_{\cal D}\right)^{2}+g\left(X_{\cal D},X_{\cal D}\right) g\left(Y_{\cal D},Y_{\cal D}\right)- g\left(\varphi X_{\cal D},Y_{\cal D}\right) ^{2}\right\}.  \label{eq-R(X,phiY,PhiY,X)}
\end{eqnarray}
Substituting $X=Z=X_{\cal D}$ and $Y=W=Y_{\cal D}$ into (\ref{eq-R(phiX,phiY)Z+R(X,Y)Z-0}), we obtain
\begin{eqnarray}
&&R\left(\varphi X_{\cal D},\varphi Y_{\cal D},X_{\cal D},Y_{\cal D}\right) - R\left(X_{\cal D},Y_{\cal D},X_{\cal D},Y_{\cal D}\right) \nonumber \\
&&\quad =\left({\alpha}^{2}-{\beta}^{2}\right) \left\{-g\left(X_{\cal D},Y_{\cal D}\right)^{2}+g\left(X_{\cal D},X_{\cal D}\right) g\left(Y_{\cal D},Y_{\cal D}\right) - g\left(\varphi X_{\cal D},Y_{\cal D}\right)^{2}\right\}. \label{eq-R(phiX,phiY)Z+R(X,Y)Z}
\end{eqnarray}
Substituting $X=X_{\cal D}+Y_{\cal D}$ in (\ref{eq-phi-sec-curv}), then using (\ref{eq-tps-Riem-cur-for-2}) and the Bianchi identity in the resulting equation, we obtain
\begin{eqnarray}
&&4R\left(X_{\cal D},\varphi X_{\cal D},X_{\cal D},\varphi Y_{\cal D}\right) +4R\left(Y_{\cal D},\varphi Y_{\cal D},Y_{\cal D},\varphi X_{\cal D}\right) +4R\left(X_{\cal D},\varphi Y_{\cal D},Y_{\cal D},\varphi X_{\cal D}\right)   \nonumber \\
&&+2R\left(\varphi X_{\cal D},\varphi Y_{\cal D},X_{\cal D},Y_{\cal D}\right) +R\left(Y_{\cal D},\varphi X_{\cal D},Y_{\cal D},\varphi X_{\cal D}\right) +R\left(X_{\cal D},\varphi Y_{\cal D},X_{\cal D},\varphi Y_{\cal D}\right)   \nonumber \\
&&\; =-c\left\{ 4g(X_{\cal D},Y_{\cal D}) ^{2}+4g(X_{\cal D},X_{\cal D}) g(X_{\cal D},Y_{\cal D}) + 4g(Y_{\cal D},Y_{\cal D}) g(X_{\cal D},Y_{\cal D})  \right. \nonumber\\
&& \qquad \quad \left. + 2g(X_{\cal D},X_{\cal D}) g(Y_{\cal D},Y_{\cal D}) \right\}.\label{eq-tps-Riem-cur-for-2a}
\end{eqnarray}
Using (\ref{eq-R(phiX,phiY)Z+R(X,Y)Z}) together with (\ref{eq-R(X,phiY,PhiY,X)}) in (\ref{eq-tps-Riem-cur-for-2a}), we get 
\begin{eqnarray}
&&4R\left(X_{\cal D},\varphi X_{\cal D},X_{\cal D},\varphi Y_{\cal D}\right) +4R\left(Y_{\cal D},\varphi Y_{\cal D},Y_{\cal D},\varphi X_{\cal D}\right) +6R\left(X_{\cal D},\varphi Y_{\cal D},Y_{\cal D},\varphi X_{\cal D}\right) +2R\left(X_{\cal D},Y_{\cal D},X_{\cal D},Y_{\cal D}\right) \nonumber\\
&&=-c\left\{4g(X_{\cal D},Y_{\cal D})^{2}+4g(X_{\cal D},X_{\cal D}) g(X_{\cal D},Y_{\cal D})+4g(Y_{\cal D},Y_{\cal D}) g(X_{\cal D},Y_{\cal D}) + 2g(X_{\cal D},X_{\cal D}) g(Y_{\cal D},Y_{\cal D}) \right\}. \label{eq-tps-Riem-cur-for-8}
\end{eqnarray}
Replacing $Y_{\cal D}$ by $-Y_{\cal D}$ in (\ref{eq-tps-Riem-cur-for-8}), we have 
\begin{eqnarray}
&&-4R\left(X_{\cal D},\varphi X_{\cal D},X_{\cal D},\varphi Y_{\cal D}\right) -4R\left(Y_{\cal D},\varphi Y_{\cal D},Y_{\cal D},\varphi X_{\cal D}\right) +6R\left(X_{\cal D},\varphi Y_{\cal D},Y_{\cal D},\varphi X_{\cal D}\right) +2R\left(X_{\cal D},Y_{\cal D},X_{\cal D},Y_{\cal D}\right)\nonumber \\
&& =-c\left\{4g(X_{\cal D},Y_{\cal D})^{2}-4g(X_{\cal D},X_{\cal D})g(X_{\cal D},Y_{\cal D}) -4g(Y_{\cal D},Y_{\cal D}) g(X_{\cal D},Y_{\cal D}) +2g(X_{\cal D},X_{\cal D}) g(Y_{\cal D},Y_{\cal D}) \right\}. \label{eq-tps-Riem-cur-for-9}
\end{eqnarray}
Adding (\ref{eq-tps-Riem-cur-for-8}) and (\ref{eq-tps-Riem-cur-for-9}), we obtain
\begin{equation}
3R\left(X_{\cal D},\varphi Y_{\cal D},Y_{\cal D},\varphi X_{\cal D}\right) + R\left(X_{\cal D},Y_{\cal D},X_{\cal D},Y_{\cal D}\right) = -c\left\{2g\left(X_{\cal D},Y_{\cal D}\right)^{2}+g\left(X_{\cal D},X_{\cal D}\right) g\left(Y_{\cal D},Y_{\cal D}\right) \right\}.\label{eq-tps-Riem-cur-for-10}
\end{equation}
Substituting $Y_{\cal D}=\varphi Y_{\cal D}$ into (\ref{eq-tps-Riem-cur-for-10}), we get 
\begin{equation}
-3R\left(X_{\cal D},Y_{\cal D},\varphi Y_{\cal D},\varphi X_{\cal D}\right) +R\left(X_{\cal D},\varphi Y_{\cal D},X_{\cal D},\varphi Y_{\cal D}\right) =-c\left\{2g\left(X_{\cal D},\varphi Y_{\cal D}\right)^{2}+g\left(X_{\cal D},X_{\cal D}\right) g\left(Y_{\cal D},Y_{\cal D}\right)\right\}.\label{eq-tps-Riem-cur-for-11}
\end{equation}
Using (\ref{eq-R(phiX,phiY)Z+R(X,Y)Z}) together with (\ref{eq-R(X,phiY,PhiY,X)}) in (\ref{eq-tps-Riem-cur-for-11}), we obtain
\begin{eqnarray}
&&3R\left(X_{\cal D},Y_{\cal D},X_{\cal D},Y_{\cal D}\right) + R\left(X_{\cal D},\varphi Y_{\cal D},Y_{\cal D},\varphi X_{\cal D}\right) = -c\left\{ 2g\left(X_{\cal D},\varphi Y_{\cal D}\right)^{2} + g\left(X_{\cal D},X_{\cal D}\right) g\left(Y_{\cal D},Y_{\cal D}\right) \right\} \nonumber \\
&&\qquad -2\left({\alpha}^{2}-{\beta}^{2}\right) \left\{-g\left(X_{\cal D},Y_{\cal D}\right)^{2}+g\left(X_{\cal D},X_{\cal D}\right) g\left(Y_{\cal D},Y_{\cal D}\right) - g\left(X_{\cal D},\varphi Y_{\cal D}\right) ^{2}\right\} .  \label{eq-tps-Riem-cur-for-12}
\end{eqnarray}
From (\ref{eq-tps-Riem-cur-for-10}) and (\ref{eq-tps-Riem-cur-for-12}), we obtain
\begin{eqnarray}
4R\left(X_{\cal D},Y_{\cal D},X_{\cal D},Y_{\cal D}\right) &=&\left(c+3\left({\alpha}^{2}-{\beta}^{2}\right) \right) \left\{g\left(X_{\cal D},Y_{\cal D}\right)^{2}-g\left(X_{\cal D},X_{\cal D}\right) g\left(Y_{\cal D},Y_{\cal D}\right) \right\}   \nonumber \\
&&-3\left(c-\left({\alpha}^{2}-{\beta}^{2}\right) \right) g\left(X_{\cal D},\varphi Y_{\cal D}\right)^{2}. \label{eq-tps-Riem-cur-for-14}
\end{eqnarray}
Substituting $Y_{\cal D}=Y_{\cal D}+W_{\cal D}$ into (\ref{eq-tps-Riem-cur-for-14}), and then using (\ref{eq-tps-Riem-cur-for-14}) in the resulting equation, we get 
\begin{eqnarray}
4R(X_{\cal D},W_{\cal D},X_{\cal D},Y_{\cal D}) &=&\left(c+3\left({\alpha}^{2}-{\beta}^{2}\right) \right) \left\{g\left(X_{\cal D},Y_{\cal D}\right) g\left(X_{\cal D},W_{\cal D}\right) -g\left(X_{\cal D},X_{\cal D}\right) g\left(Y_{\cal D},W_{\cal D}\right) \right\}  \nonumber \\
&&-3\left(c-\left({\alpha}^{2}-{\beta}^{2}\right) \right) g\left(X_{\cal D},\varphi Y_{\cal D}\right) g\left(X_{\cal D},\varphi W_{\cal D}\right).  \label{eq-a}
\end{eqnarray}
Substituting $X_{\cal D}=X_{\cal D}+Z_{\cal D}$ into (\ref{eq-tps-Riem-cur-for-14}), and then using (\ref{eq-tps-Riem-cur-for-14}) in the resulting equation, we get
\begin{eqnarray}
4R\left(X_{\cal D},Y_{\cal D},Z_{\cal D},Y_{\cal D}\right) &=& \left(c+3\left({\alpha}^{2}-{\beta}^{2}\right) \right) \left\{ g\left(X_{\cal D},Y_{\cal D}\right) g\left(Y_{\cal D},Z_{\cal D}\right) - g\left(Y_{\cal D},Y_{\cal D}\right)g\left(X_{\cal D},Z_{\cal D}\right) \right\}   \nonumber \\
&&-3\left(c-\left({\alpha}^{2}-{\beta}^{2}\right) \right) g\left(\varphi Y_{\cal D},Z_{\cal D}\right) g\left(X_{\cal D},\varphi Y_{\cal D}\right) .  \label{eq-b}
\end{eqnarray}
Substituting $X_{\cal D}=X_{\cal D}+Z_{\cal D}$ and $Y_{\cal D}=Y_{\cal D}+W_{\cal D}$ into (\ref{eq-tps-Riem-cur-for-14}), then using (\ref{eq-tps-Riem-cur-for-14}), (\ref{eq-a}) and (\ref{eq-b}) in the resulting equation, we get 
\begin{eqnarray}
&&4R\left(X_{\cal D},Y_{\cal D},Z_{\cal D},W_{\cal D}\right) + 4R\left(X_{\cal D},W_{\cal D},Z_{\cal D},Y_{\cal D}\right)\nonumber \\
&&\quad =\left(c+3\left({\alpha}^{2}-{\beta}^{2}\right) \right)\left\{g\left(Y_{\cal D},Z_{\cal D}\right) g\left(X_{\cal D},W_{\cal D}\right) + g\left(X_{\cal D},Y_{\cal D}\right)g\left(Z_{\cal D},W_{\cal D}\right) - 2g\left(X_{\cal D},Z_{\cal D}\right) g\left(Y_{\cal D},W_{\cal D}\right) \right\}   \nonumber \\
&&\qquad -3\left(c-\left({\alpha}^{2}-{\beta}^{2}\right) \right) \left\{g\left(\varphi Y_{\cal D},Z_{\cal D}\right) g\left(X_{\cal D},\varphi W_{\cal D}\right) +g\left(X_{\cal D},\varphi Y_{\cal D}\right) g\left(Z_{\cal D},\varphi W_{\cal D}\right) \right\}. \label{eq-tps-Riem-cur-for-15}
\end{eqnarray}
Interchanging $Y_{\cal D}$ with $Z_{\cal D}$ in (\ref{eq-tps-Riem-cur-for-15}), we get
\begin{eqnarray}
&&4R\left(X_{\cal D},Z_{\cal D},Y_{\cal D},W_{\cal D}\right) + 4R\left(X_{\cal D},W_{\cal D},Y_{\cal D},Z_{\cal D}\right)  \nonumber \\
&&\quad =\left(c+3\left({\alpha}^{2}-{\beta}^{2}\right) \right) \left\{g\left(Y_{\cal D},Z_{\cal D}\right) g\left(X_{\cal D},W_{\cal D}\right) +g\left(X_{\cal D},Z_{\cal D}\right) g\left(Y_{\cal D},W_{\cal D}\right) -2g\left(X_{\cal D},Y_{\cal D}\right) g\left(Z_{\cal D},W_{\cal D}\right) \right\}   \nonumber \\
&&\qquad -3\left(c-\left({\alpha}^{2}-{\beta}^{2}\right) \right) \left\{g\left(Y_{\cal D},\varphi Z_{\cal D}\right) g\left(X_{\cal D},\varphi W_{\cal D}\right) +g\left(X_{\cal D},\varphi Z_{\cal D}\right) g\left(Y_{\cal D},\varphi W_{\cal D}\right) \right\}. \label{eq-tps-Riem-cur-for-16}
\end{eqnarray}
Subtracting (\ref{eq-tps-Riem-cur-for-16}) from (\ref{eq-tps-Riem-cur-for-15}) and using Bianchi identity, and interchanging $Y_{\cal D}$ with $W_{\cal D}$ in the resulting equation, we get
\begin{eqnarray}
&&4R\left(X_{\cal D},Y_{\cal D},Z_{\cal D},W_{\cal D}\right) 
=\left(c+3\left({\alpha}^{2}-{\beta}^{2}\right) \right) \left\{g\left(Y_{\cal D},Z_{\cal D}\right) 
g\left(X_{\cal D},W_{\cal D}\right) -g\left(X_{\cal D},Z_{\cal D}\right) g\left(Y_{\cal D},W_{\cal D}\right) 
\right\}   \nonumber \\
&&\quad \qquad \qquad\qquad\qquad \qquad -\left(c-\left({\alpha}^{2}-{\beta}^{2}\right) \right) \left\{g\left(\varphi 
Y_{\cal D},Z_{\cal D}\right) g\left(X_{\cal D},\varphi W_{\cal D}\right) -g\left(X_{\cal D},\varphi 
Z_{\cal D}\right) g\left(\varphi Y_{\cal D},W_{\cal D}\right) \right.  \nonumber\\
&&\qquad \qquad \qquad \qquad \qquad \qquad \left. +2g\left(X_{\cal D},\varphi Y_{\cal D}\right) 
g\left(Z_{\cal D},\varphi W_{\cal D}\right) \right\} .  \label{eq-tps-Riem-cur-for-17a}
\end{eqnarray}
Substituting $X_{\cal D}=X-\eta \left(X\right) \xi $,\quad $Y_{\cal D}=Y-\eta \left(Y\right) \xi $, \quad $Z_{\cal D}=Z-\eta \left(Z\right)\xi $,\quad $W_{\cal D}=W-\eta \left(W\right) \xi $ into (\ref{eq-tps-Riem-cur-for-17a}), we get
\begin{eqnarray}
&&R\left(X,Y,Z,W\right) -\eta \left(W\right) R\left(X,Y,Z,\xi \right) -\eta \left(Z\right) R\left(X,Y,\xi ,W\right) -\eta \left(Y\right) R\left(X,\xi,Z,W\right)   \nonumber \\
&&+\eta \left(Y\right) \eta \left(W\right) R\left(X,\xi,Z,\xi \right)+\eta \left(Z\right) \eta \left(Y\right) R\left(X,\xi ,\xi ,W\right)-\eta \left(X\right) R\left(\xi ,Y,Z,W\right)   \nonumber \\
&&+\eta \left(X\right) \eta \left(W\right) R\left(\xi,Y,Z,\xi \right) +\eta \left(X\right) \eta \left(Z\right) R\left(\xi,Y,\xi,W\right) \nonumber \\
&&\quad =\frac{1}{4}\left(c+3\left({\alpha}^{2}-{\beta}^{2}\right) \right) \left\{ g\left(Y,Z\right) g\left(X,W\right) -g\left(X,Z\right) g\left(Y,W\right) \right.   \nonumber \\
&&\qquad \qquad \qquad \qquad \quad \qquad \left. -g\left(Y,Z\right) \eta \left(X\right) \eta \left(W\right) +g\left(X,Z\right) \eta \left(Y\right) \eta \left(W\right) \right.   \nonumber \\
&&\qquad \qquad \qquad \qquad \qquad \quad \left. -\eta \left(Y\right) \eta \left(Z\right) g\left(X,W\right) +\eta \left(X\right) \eta \left(Z\right) g\left(Y,W\right) \right\}   \nonumber \\
&&\qquad -\frac{1}{4}\left(c-\left({\alpha}^{2}-{\beta}^{2}\right) \right) 
\left\{-g\left(X,\varphi{W}\right) g\left(\varphi{Z},Y\right) +g\left(X,\varphi{Z}\right) 
g\left(\varphi{W},Y\right) \right.   \nonumber \\
&&\qquad \qquad \qquad \qquad \qquad \quad \left. +2g\left(X,\varphi {Y}\right) 
g\left(Z,\varphi{W}\right) \right\} . \label{eq-tps-Riem-cur-for-19}
\end{eqnarray}
Using (\ref{eq-tr-sasa-cur-ope}) and (\ref{eq-funda-2-form}) in (\ref{eq-tps-Riem-cur-for-19}), and taking (\ref{eq-V-Sym-Product}) and (\ref{eq-V-KN-Product}) into account, we get
\begin{eqnarray}
R &=&\frac{c+3\left({\alpha}^{2}-{\beta}^{2}\right) }{8}g\circledast g \nonumber \\
&&+\frac{c-\left({\alpha}^{2}-{\beta}^{2}\right) }{4}\left(\frac{1}{2}{\varphi}^{\flat} \circledast {\varphi}^{\flat} -{\varphi}^{\flat} \circledcirc {\varphi}^{\flat} -g\circledast \left(\eta \otimes \eta \right) \right)   \nonumber \\
&&+{\rm d}\beta \left(\xi \right) g\circledast \left(\eta \otimes \eta \right) -g\circledast \left({\rm d}\beta \otimes \eta +\eta \otimes {\rm d}\beta \right)   \nonumber \\
&&-{\varphi}^{\flat} \circledast \left({\rm d}\alpha \otimes \eta -\eta \otimes {\rm d}\alpha \right).  \label{eq-tps-Riem-cur-for-30}
\end{eqnarray}
If ${\rm dim}(M)=3$, then using (\ref{eq-3D-tsm-R(Y,Y)xi}) in (\ref{eq-tps-Riem-cur-for-19}), we get 
\begin{eqnarray}
R &=&\frac{c+3\left({\alpha}^{2}-{\beta}^{2}\right) }{8}g\circledast g \nonumber \\
&&+\frac{c-\left({\alpha}^{2}-{\beta}^{2}\right) }{4}\left(\frac{1}{2}{\varphi}^{\flat} \circledast {\varphi}^{\flat} -{\varphi}^{\flat} \circledcirc {\varphi}^{\flat} -g\circledast \left(\eta \otimes \eta \right) \right)   \nonumber \\
&&+{\rm d}\beta \left(\xi \right) g\circledast \left(\eta \otimes \eta \right) -g\circledast \left({\rm d}\beta \otimes \eta +\eta \otimes {\rm d}\beta \right)   \nonumber \\
&&+g\circledast \left(\left(\varphi {\rm grad}\left(\alpha \right) \right)^{\sharp}\otimes \eta +\eta \otimes \left(\varphi {\rm grad}\left(\alpha \right) \right)^{\sharp}\right) .  \label{eq-tps-Riem-cur-for-31}
\end{eqnarray}
For ${\rm dim}(M)=3$, comparing (\ref{eq-tps-Riem-cur-for-30}) and (\ref{eq-tps-Riem-cur-for-31}), we obtain
\begin{equation}
{\varphi}^{\flat} \circledast \left({\rm d}\alpha \otimes \eta -\eta \otimes {\rm d}\alpha \right) =-g\circledast \left(\left(\varphi {\rm grad}\left(\alpha\right) \right)^{\sharp}\otimes \eta +\eta \otimes \left(\varphi {\rm grad}\left(\alpha \right) \right)^{\sharp}\right). \label{eq-tps-Riem-cur-for-32}
\end{equation}
Next, if ${\rm dim}(M)>3$, then $\alpha$ reduces to be a constant, therefore (\ref{eq-tps-Riem-cur-for-32}) holds identically. Using (\ref{eq-tps-Riem-cur-for-32}) in (\ref{eq-tps-Riem-cur-for-30}), we get (\ref{eq-tps-Riem-cur-for-33}).
Conversely, for a unit tangent vector field $X$ orthogonal to $\xi$, equation (\ref{eq-tps-Riem-cur-for-33}), yields
\[
R\left(X,\varphi {X},\varphi {X},X\right) =c, 
\]
which establishes the stated result. $\blacksquare$

As an immediate consequence of the curvature formula \eqref{eq-tps-Riem-cur-for-33}, we obtain the corresponding curvature formulas for the important subclasses of trans-Sasakian manifolds introduced in Remark~\ref{rem-subcls-tsm}.

\begin{cor}
In a smooth manifold $M$, let $(\varphi,\xi,\eta,g)$ be a ${\cal C}_{6}$-structure. For $M$, the necessary and sufficient condition of having pointwise constant $\varphi$-sectional curvature $c$ is that its Riemannian-Christoffel curvature tensor $R$ is given by  
\begin{eqnarray}
R &=&\frac{c+3\alpha ^{2}}{8}g\circledast g+\frac{c-\alpha ^{2}}{4}\left( 
\frac{1}{2}{\varphi}^{\flat} \circledast {\varphi}^{\flat} -{\varphi}^{\flat} \circledcirc {\varphi}^{\flat} -g\circledast
\left( \eta \otimes \eta \right) \right)  \nonumber \\
&&+g\circledast \left( \left( \varphi {\rm grad}\left( \alpha \right)
\right) ^{\sharp }\otimes \eta +\eta \otimes \left( \varphi {\rm grad}\left(
\alpha \right) \right) ^{\sharp }\right).  \label{eq-C6-sf}
\end{eqnarray}
\end{cor}

\begin{cor}
{\rm {\cite[Janssens and Vanhecke 1981]{Janssens-Vanhecke-81-KMJ}}} In a smooth manifold $M$, let $(\varphi,\xi,\eta,g)$ be an $\alpha$-Sasakian structure. For $M$, the necessary and sufficient condition of having pointwise constant $\varphi$-sectional curvature $c$ is that its Riemannian-Christoffel curvature tensor $R$ is given by
\begin{equation}
R=\frac{c+3\alpha ^{2}}{8}g\circledast g+\frac{c-\alpha ^{2}}{4}\left(\frac{1%
}{2}{\varphi}^{\flat} \circledast {\varphi}^{\flat} -{\varphi}^{\flat} \circledcirc {\varphi}^{\flat} -g\circledast \left(\eta
\otimes \eta \right) \right) .  \label{eq-alpha-Sasakian-sf}
\end{equation}
\end{cor}

\begin{cor}
{\rm {(\cite[Blair 2010]{Blair-10-Book}, \cite[Ogiue 1964]{Ogiue-64-Kodai})}} In a smooth manifold $M$, let $(\varphi,\xi,\eta,g)$ be an Sasakian structure. For $M$, the necessary and sufficient condition of having pointwise constant $\varphi$-sectional curvature $c$ is that its Riemannian-Christoffel curvature tensor $R$ is given by
\begin{equation}
R=\frac{c+3}{8}g\circledast g+\frac{c-1}{4}\left(\frac{1}{2}{\varphi}^{\flat} \circledast
{\varphi}^{\flat} -{\varphi}^{\flat} \circledcirc {\varphi}^{\flat} -g\circledast \left(\eta \otimes \eta \right)
\right).  \label{eq-Sasakian-sf}
\end{equation}
\end{cor}

\begin{cor}
{\rm (cf.~\cite[Wu et al. 2024, Theorem 3.1]{Wu-Sheng-Zhang-24-MJM})} In a smooth manifold $M$, let $(\varphi,\xi,\eta,g)$ be an ${\cal C}_{5}$-structure. For $M$, the necessary and sufficient condition of having pointwise constant $\varphi$-sectional curvature $c$ is that its Riemannian-Christoffel curvature tensor $R$ is given by
\begin{eqnarray}
R &=&\frac{c-3\beta ^{2}}{8}g\circledast g+\frac{c+\beta ^{2}}{4}\left(\frac{1}{2}{\varphi}^{\flat} \circledast {\varphi}^{\flat} -{\varphi}^{\flat} \circledcirc {\varphi}^{\flat} -g\circledast\left(\eta \otimes \eta \right) \right)  \nonumber \\
&&+d\beta \left(\xi \right) g\circledast \left(\eta \otimes \eta \right)
-g\circledast \left(d\beta \otimes \eta +\eta \otimes d\beta \right).
\label{eq-f-Kenmotsu-sf}
\end{eqnarray}
\end{cor}

\begin{cor}
{\rm {\cite[Janssens and Vanhecke 1981]{Janssens-Vanhecke-81-KMJ}}} In a smooth manifold $M$, let $(\varphi,\xi,\eta,g)$ be an $\beta$-Kenmotsu structure. For $M$, the necessary and sufficient condition of having pointwise constant $\varphi$-sectional curvature $c$ is that its Riemannian-Christoffel curvature tensor $R$ is given by
\begin{equation}
R=\frac{c-3\beta ^{2}}{8}g\circledast g+\frac{c+\beta ^{2}}{4}\left(\frac{1}{%
2}{\varphi}^{\flat} \circledast {\varphi}^{\flat} -{\varphi}^{\flat} \circledcirc {\varphi}^{\flat} -g\circledast \left(\eta
\otimes \eta \right) \right).  \label{eq-beta-Kenmotsu-sf}
\end{equation}
\end{cor}

\begin{cor}
{\rm \cite[Kenmotsu 1972]{Kenmotsu-72-Tohoku}} In a smooth manifold $M$, let $(\varphi,\xi,\eta,g)$ be an Kenmotsu structure. For $M$, the necessary and sufficient condition of having pointwise constant $\varphi$-sectional curvature $c$ is that its Riemannian-Christoffel curvature tensor $R$ is given by
\begin{equation}
R=\frac{c-3}{8}g\circledast g+\frac{c+1}{4}\left(\frac{1}{2}{\varphi}^{\flat} \circledast
{\varphi}^{\flat} -{\varphi}^{\flat} \circledcirc {\varphi}^{\flat} +g\circledast (\eta \otimes \eta )\right) .
\label{eq-Kenmotsu-sf}
\end{equation}
\end{cor}

\begin{cor}
{\rm {\cite[Ludden 1970]{Ludden-70-JDG}}} In a smooth manifold $M$, let $(\varphi,\xi,\eta,g)$ be an cosymplectic structure. For $M$, the necessary and sufficient condition of having pointwise constant $\varphi$-sectional curvature $c$ is that its Riemannian-Christoffel curvature tensor $R$ is given by
\begin{equation}
R=\frac{c}{8}g\circledast g+\frac{c}{4}\left(\frac{1}{2}{\varphi}^{\flat} \circledast
{\varphi}^{\flat} -{\varphi}^{\flat} \circledcirc {\varphi}^{\flat} -g\circledast \left(\eta \otimes \eta \right)
\right) .  \label{eq-cosymplectic-sf}
\end{equation}
\end{cor}

Motivated by the curvature formula \eqref{eq-tps-Riem-cur-for-33}, we introduce the following definition.
\begin{defn-new}
A trans-Sasakian manifold $(M,\varphi,\xi,\eta,g)$ of type $(\alpha,\beta)$ is called a trans-Sasakian space form if it has pointwise constant $\varphi$-sectional curvature $c$. Such a manifold is denoted by $M(c)$.
\end{defn-new}

\begin{rem-new}
A trans-Sasakian space form $M(c)$ is an example of a Blair space form \cite[Tripathi 2025, Eq. (4.2.4)]{Tri-25-Preprint} with 
\[
T^{\left(g\right)}=g\circledast \left(\left(\varphi {\rm grad}\left(\alpha \right) \right)^{\sharp}\otimes \eta +\eta \otimes \left(\varphi {\rm grad}\left(\alpha \right) \right)^{\sharp}-\left({\rm d}\beta \otimes \eta +\eta \otimes {\rm d}\beta \right) \right),
\]
\[
c^{\left(g\right)}=\frac{c+3\left({\alpha}^{2}-{\beta}^{2}\right) }{4},\qquad c^{\left({\varphi}^{\flat} \right) }=\frac{c-\left({\alpha}^{2}-\beta^{2}\right) }{4},\qquad c^{\left(g,\eta \right) }={\rm d}\beta \left(\xi\right) -\frac{c-\left({\alpha}^{2}-{\beta}^{2}\right)}{4}. 
\]
\end{rem-new}

\begin{prop}
In a trans-Sasakian space form $M(c)$, the Ricci tensor $S$ and scalar curvature $\tau$, defined by {\rm(\ref{eq-scal-cur})}, are expressed as follows: 
\begin{eqnarray}
S\left(X,Y\right)  &=&\left(\frac{\left(c+3\left({\alpha}^{2}-\beta^{2}\right) \right)n +\left(c-\left({\alpha}^{2}-{\beta}^{2}\right) \right) }{2}-{\rm d}\beta \left(\xi \right) \right) g\left(X,Y\right)   \nonumber \\
&&-\left(\frac{\left(n+1\right) \left(c-\left({\alpha}^{2}-\beta^{2}\right) \right) }{2}-\left(2n-1\right) {\rm d}\beta(\xi)\right) \eta(X)\eta(Y)   \nonumber \\
&&-\left(2n-1\right) \left\{ \eta \left(X\right) {\rm d}\beta \left(Y\right) +{\rm d}\beta \left(X\right) \eta \left(Y\right) \right. \nonumber \\
&&\qquad \qquad \quad \left. +\eta \left(X\right) {\rm d}\alpha \left(\varphi{Y}\right) +{\rm 
d}\alpha \left(\varphi{X}\right) \eta \left(Y\right)\right\},  \label{eq-Ricci-tensor-tsf}
\end{eqnarray}
\begin{equation}
\tau =\frac{1}{4}\left\{\left(2n^{2}+n\right) \left(c+3\left(\alpha^{2}-{\beta}^{2}\right) \right) +\left(c-\left({\alpha}^{2}-\beta^{2}\right) \right)n \right\} -2n{\rm d}\beta \left(\xi \right) \label{eq-scalar-curvature-tsf}
\end{equation}
for all vector fields $X$, $Y$ on $M(c)$.
\end{prop}

\noindent {\bf Proof}. By a straightforward calculation from (\ref{eq-tps-Riem-cur-for-33}), we obtain (\ref{eq-Ricci-tensor-tsf}) and (\ref{eq-scalar-curvature-tsf}). $\blacksquare$

\begin{thm} \label{th-sec-3D-tssf}
Every $3$-dimensional trans-Sasakian manifold $\left(M,\varphi,\xi,\eta,g\right)$ has pointwise constant $\varphi$-sectional curvature $c$, given by 
\begin{equation}
c=\tau-2\left({\alpha}^{2}-{\beta}^{2}-{\rm d}\beta \left(\xi \right)\right),  \label{eq-3D-tsm-c}
\end{equation}
where $\tau$ is the scalar curvature of $M$ defined by {\rm (\ref{eq-scal-cur})}.
\end{thm}

\noindent {\bf Proof.} By taking inner product of (\ref{eq-3D-tsm-R(X,Y,Z)}) with $W$, then substituting $Y=Z=\varphi{X}$ and $W=X$ into the resulting equation, we obtain 
\[
c=R\left(X,\varphi{X},\varphi{X},X\right) =\tau -2\left({\alpha}^{2}-\beta^{2}-{\rm d}\beta 
\left(\xi \right) \right),
\]
which proves (\ref{eq-3D-tsm-c}). $\blacksquare$

\begin{ex-new}\label{ex-tsm-expo}
Consider a $3$-dimensional smooth manifold $M={\Bbb R}^{3}$ with natural coordinates $\left(x^{1},x^{2},x^{3}\right)$ in ${\Bbb R}^{3}$. Let $\{U_{1},U_{2},U_{3}\}$ be a global frame $\{U_{1},U_{2},U_{3}\}$ on $M$, defined by
\[
U_{1}=e^{x^{3}}\left( \frac{\partial }{\partial x^{1}}+x^{2}\frac{\partial }{\partial x^{3}}\right) ,\quad U_{2}=e^{x^{3}}\frac{\partial }{\partial x^{2}},\quad U_{3}=\frac{\partial }{\partial x^{3}}. 
\]
Note that $\{U_{1},U_{2},U_{3}\}$ is linearly independent. Suppose that $g$ is a Riemannian metric, given by
\begin{eqnarray*}
g(U_{i},U_{j}) = \delta_{ij}, \quad i,j \in \{1,2,3\},
\end{eqnarray*}
where $\delta_{ij}$ is the well-known Kronecker delta.
Suppose that $\eta$ is the $1$-form on $M$, defined by
\[
\eta (X)=g(X,U_{3})
\]
for every vector field $X$ on $M$. Next, suppose that $\varphi$ is a $(1,1)$ type tensor field, defined by
\[
\varphi (U_{1})=U_{2},\quad \varphi (U_{2})=-U_{1},\quad \varphi (U_{3})=0. 
\]
Let $\nabla$ denote the Levi-Civita connection of $g$, and $R$ be the Riemann--Christoffel curvature tensor of $M$. Then, we have
\[
[U_{1},U_{2}]=x^{2}e^{x^{3}}U_{2}-e^{2x^{3}}U_{3},\quad [U_{2},U_{3}]=-U_{2},\quad [U_{1},U_{3}]=-U_{1}.  
\]
Also, we have
\begin{equation*}
\begin{aligned}
\nabla_{U_{1}}U_{1} &= U_{3}, 
&\qquad
\nabla_{U_{1}}U_{2} &= -\frac{1}{2}e^{2x^{3}}U_{3}, 
&\qquad
\nabla_{U_{1}}U_{3} &= -U_{1}+\frac{1}{2}e^{2x^{3}}U_{2}, \\[0.5em]
\nabla_{U_{2}}U_{1} &= \frac{1}{2}e^{2x^{3}}U_{3}-x^{2}e^{x^{3}}U_{2}, 
&\qquad
\nabla_{U_{2}}U_{2} &= U_{3}+x^{2}e^{x^{3}}U_{1}, 
&\qquad
\nabla_{U_{2}}U_{3} &= -U_{2}-\frac{1}{2}e^{2x^{3}}U_{1}, \\[0.5em]
\nabla_{U_{3}}U_{1} &= \frac{1}{2}e^{2x^{3}}U_{2}, 
&\qquad
\nabla_{U_{3}}U_{2} &= -\frac{1}{2}e^{2x^{3}}U_{1}, 
&\qquad
\nabla_{U_{3}}U_{3} &= 0.
\end{aligned}
\end{equation*}
Then, $(\varphi,\xi,\eta,g)$ defines a trans Sasakian structure on $M$
with $\left(\alpha,\beta\right) =\left(-\frac{1}{2}e^{2x^{3}},-1\right)$ \cite[Shaikh {\em et al.} 2008]{Shaikh-Baishya-Eyasmin-08-DM}. Making use of the previous expressions, we arrive at the following components of $R$.
\begin{equation*}
\begin{aligned}
R(U_{1},U_{2})U_{1} &= \left( \frac{3}{4}e^{4x^{3}}+1\right)U_{2}
+x^{2}e^{3x^{3}}U_{3}, 
&\qquad
R(U_{1},U_{2})U_{2} &= -\left( \frac{3}{4}e^{4x^{3}}+1\right)U_{1}, \\[0.5em]
R(U_{1},U_{2})U_{3} &= -x^{2}e^{3x^{3}}U_{1},
&\qquad
R(U_{1},U_{3})U_{1} &= x^{2}e^{3x^{3}}U_{2}
+\left( 1-\frac{1}{4}e^{4x^{3}}\right)U_{3}, \\[0.5em]
R(U_{1},U_{3})U_{2} &= -x^{2}e^{3x^{3}}U_{1},
&\qquad
R(U_{1},U_{3})U_{3} &= \left( \frac{1}{4}e^{4x^{3}}-1\right)U_{1}, \\[0.5em]
R(U_{2},U_{3})U_{1} &= 0,
&\qquad
R(U_{2},U_{3})U_{2} &= \left( 1-\frac{1}{4}e^{4x^{3}}\right)U_{3}, \\[0.5em]
R(U_{2},U_{3})U_{3} &= \left( \frac{1}{4}e^{4x^{3}}-1\right)U_{2}.
\end{aligned}
\end{equation*}
A straightforward calculation based on the curvature tensor gives the following entries of the Ricci tensor $S$.
\begin{eqnarray*}
S\left( U_{1},U_{1}\right) &=&-\frac{1}{2}e^{4x^{3}}-2, \qquad
S\left( U_{2},U_{2}\right) =-\frac{1}{2}e^{4x^{3}}-2, \\
S\left( U_{3},U_{3}\right) &=&\frac{1}{2}e^{4x^{3}}-2.
\end{eqnarray*}
Using the above expressions, we find that the scalar curvature $\tau$ as follows:
\[
\tau = -\frac{1}{4}e^{4x^{3}}-3. 
\]
A direct computations from the curvature tensor yields the $\varphi$-sectional curvature of $M$ as
\[
K\left(U_1,\varphi {U_1}\right)= K\left(U_2,\varphi {U_2}\right) = -\frac{3}{4}e^{4x^{3}}-1.
\]
From the above computations, it follows that $M$ satisfies relations (\ref{eq-tps-Riem-cur-for-33}) and  (\ref{eq-3D-tsm-c}). Consequently, Theorems \ref{th-ts-sf} and \ref{th-sec-3D-tssf} are verified. 
\end{ex-new}

\section{Trans-Sasakian statistical manifolds\label{sec-trans-sasa-stat}}
In this section, we present several results concerning trans-Sasakian statistical manifolds and introduce the notion of trans-Sasakian statistical space forms.

Let $(M,\varphi,\xi,\eta,g)$ be a $(2n+1)$-dimensional almost contact metric manifolds and let $({\sc D},g)$ defines statistical structure on $M$. If
\begin{equation}
{\sc B}^{\sc d}_{X}\varphi Y +\varphi {\sc B}^{\sc d}_{X}Y = 0
\label{eq-almost-con-statis}
\end{equation}
for all vector fields $X$, $Y$ on $M$, then $(M,\varphi,\xi,\eta,g,{\sc D})$ is called almost contact metric statistical manifold (see \cite[Akbari and Malek 2021, Definition~3]{Akbari-Malek-21-DGA}).

\begin{defn-new}\cite[Kazan 2019, Definition~1]{Kazan-19-PhyA}
An almost contact metric statistical structure on a trans-Sasakian manifold of type $(\alpha,\beta)$ is called a $(\alpha,\beta)$ trans-Sasakian statistical structure, and such a manifold is called a trans-Sasakian statistical manifold of type $(\alpha,\beta)$.
\end{defn-new}

\begin{rem-new} \label{rem-class-tssm}
An $(\alpha,\beta)$ trans-Sasakian statistical structure reduces to
\begin{itemize}
\item[{\bf (a)}] ${\cal C}_{6}$-statistical structure if $\beta=0$,

\item[{\bf (b)}] an $\alpha$-Sasakian statistical structure if $\beta=0\neq \alpha =$ constant,
    
\item[{\bf (c)}] a Sasakian statistical structure \cite[Furuhata {\em et al.} 2017]{Furuhata-Hasegawa-Okuyama-Sato-Shahid-17-JGP} if $\beta=0$, $\alpha=1$,

\item[{\bf (d)}] ${\cal C}_{5}$-statistical structure (cf.~\cite[Wu {\em et al.} 2024]{Wu-Sheng-Zhang-24-MJM}) if $\alpha=0$,

\item[{\bf (e)}] a $\beta$-Kenmotsu statistical structure if $\alpha =0\neq \beta =$ constant,

\item[{\bf (f)}] a Kenmotsu statistical structure \cite[Furuhata {\em et al.} 2017]{Furu-Hasegawa-Okuyama-Sato-17-JG} if $\alpha=0$, $\beta=1$,

\item[{\bf (g)}] a cosymplectic statistical structure \cite[Malek and Akbari 2020]{Malek-Akbari-20-BIMS} if $\alpha =0=\beta$.
\end{itemize}
\end{rem-new}

%
%

%

\begin{prop}
Let $(M,\varphi,\xi,\eta,g,{\sc D})$ be an almost contact metric statistical manifold. Then
\begin{equation}
\left({\sc D}_{X}\varphi \right) Y-\left({\sc D}_{Y}\varphi \right)
X=\left(\nabla_{X}\varphi \right) Y-\left(\nabla_{Y}\varphi
\right) X,  \label{eq-nabla(X,varphiY)=nabla0(X,varphiY)}
\end{equation}
\begin{equation}
\left({\sc D}_{X}^{\ast}\varphi \right) Y-\left({\sc D}_{Y}^{\ast}\varphi \right) X=\left(\nabla_{X}\varphi \right) Y-\left(\nabla_{Y}\varphi \right) X,
\label{eq-nabla(X,varphiY)=nabla0(X,varphiY)-star}
\end{equation}
\begin{equation}
\left({\sc D}_{X}\eta \right) Y=g\left(Y,{\sc D}_{X}^{\ast}\xi \right) ,
\label{eq-nablaXeta-2}
\end{equation}
\begin{equation}
\left({\sc D}_{X}^{\ast}\eta \right) Y=g\left(Y,{\sc D}_{X}\xi \right)
\label{eq-nablaXeta-star}
\end{equation}
for all vector fields $X$, $Y$ on $M$.
\end{prop}

\noindent {\bf Proof.} From (\ref{eq-K(X,Y)}), we have
\begin{eqnarray}
{\sc B}^{\sc d}_{X}Y &=&{\sc D}_{X}Y-\nabla_{X}Y,  \label{eq-K(X,Y)-1} \\
{\sc B}^{\sc d}_{Y}X &=&{\sc D}_{Y}X-\nabla_{Y}X  \label{eq-K(Y,X)-1}
\end{eqnarray}
for all vector fields $X$, $Y$ on $M$. Subtracting (\ref{eq-K(Y,X)-1}) from (\ref{eq-K(X,Y)-1}), and using (\ref{eq-K-properties}), we get
\begin{equation}
{\sc D}_{X}Y-{\sc D}_{Y}X=\nabla_{X}Y-\nabla_{Y}X.
\label{eq-nabla(X,Y)=nabla0(X,Y)}
\end{equation}
Replacing $Y$ by $\varphi{Y}$ in (\ref{eq-K(X,Y)-1}) and $X$ by $\varphi{X}$ in (\ref{eq-K(Y,X)-1}), we obtain
\begin{equation}
{\sc B}^{\sc d}_{X}\varphi Y = {\sc D}_{X}\varphi Y-\nabla_{X}\varphi Y,
\label{eq-K(X,Y)-2} 
\end{equation}
\begin{equation}
{\sc B}^{\sc d}_{Y}\varphi X = {\sc D}_{Y}\varphi X-\nabla_{Y}\varphi X.
\label{eq-K(Y,X)-2}
\end{equation}
Subtracting (\ref{eq-K(Y,X)-2}) from (\ref{eq-K(X,Y)-2}), and using (\ref{eq-K-properties}) and (\ref{eq-almost-con-statis}), we get
\begin{equation}
{\sc D}_{X}\varphi Y-{\sc D}_{Y}\varphi X=\nabla_{X}\varphi Y-\nabla_{Y}\varphi X.  \label{eq-nabla(X,varphiY)=nabla0(X,varphiY)-1}
\end{equation}
From (\ref{eq-nabla(X,Y)=nabla0(X,Y)}) and (\ref{eq-nabla(X,varphiY)=nabla0(X,varphiY)-1}), we get (\ref{eq-nabla(X,varphiY)=nabla0(X,varphiY)}).
From (\ref{eq-nabla-star-zero}), we have
\begin{equation}
{\sc D}_{X}Y+{\sc D}_{X}^{\ast}Y=2\nabla_{X}Y.
\label{eq-nablaXY+star-1}
\end{equation}
\begin{equation}
{\sc D}_{Y}X+{\sc D}_{Y}^{\ast}X=2\nabla_{Y}X.
\label{eq-nablaXY+star-12}
\end{equation}
Subtracting (\ref{eq-nablaXY+star-12}) from (\ref{eq-nablaXY+star-1}), and using (\ref{eq-nabla(X,Y)=nabla0(X,Y)}), we get
\begin{equation}
{\sc D}_{X}^{\ast}Y - {\sc D}_{Y}^{\ast}X = \nabla_{X}Y-\nabla_{Y}X.
\label{eq-nablaXY+star-3}
\end{equation}
Replacing $Y$ by $\varphi Y$ in (\ref{eq-nablaXY+star-1}) and $X$ by $\varphi X$ in (\ref{eq-nablaXY+star-12}), we obtain
\begin{equation}
{\sc D}_{X}\varphi Y + {\sc D}_{X}^{\ast}\varphi Y=2\nabla_{X}\varphi Y,
\label{eq-nablaXY+star-4}
\end{equation}
\begin{equation}
{\sc D}_{Y}\varphi X + {\sc D}_{Y}^{\ast}\varphi X=2\nabla_{Y}\varphi X.
\label{eq-nablaXY+star-5}
\end{equation}
Subtracting (\ref{eq-nablaXY+star-5}) from (\ref{eq-nablaXY+star-4}), and using (\ref{eq-nabla(X,varphiY)=nabla0(X,varphiY)-1}), we get
\begin{equation}
{\sc D}_{X}^{\ast}\varphi Y - {\sc D}_{Y}^{\ast}\varphi X=\nabla_{X}\varphi Y-\nabla_{Y}\varphi X.  \label{eq-nablaXY+star-6}
\end{equation}
From (\ref{eq-nablaXY+star-3}) and (\ref{eq-nablaXY+star-6}), we obtain (\ref{eq-nabla(X,varphiY)=nabla0(X,varphiY)-star}).\\
Now, we have
\begin{equation}
\left({\sc D}_{X}\eta \right) Y = X\eta (Y)-\eta \left({\sc D}_{X}Y\right).
\label{eq-nablaXeta-0}
\end{equation}
Putting $Z=\xi$ in (\ref{eq-nabla-star}), we get
\begin{equation}
X\eta \left(Y\right) =g\left({\sc D}_{X}Y,\xi \right) + g\left(Y,{\sc D}_{X}^{\ast}\xi \right).  \label{eq-nablaXeta-1}
\end{equation}
Using (\ref{eq-nablaXeta-1}) in (\ref{eq-nablaXeta-0}), we get (\ref{eq-nablaXeta-2}). Next, we have
\begin{equation}
\left({\sc D}_{X}^{\ast}\eta \right) Z = X\eta (Z)-\eta \left({\sc D}_{X}^{\ast}Z\right).  \label{eq-nablaXeta-star-1}
\end{equation}
Putting $Y=\xi$ in (\ref{eq-nabla-star}), we get
\begin{equation}
X\eta (Z)=g\left({\sc D}_{X}\xi ,Z\right) + g\left(\xi,{\sc D}_{X}^{\ast}Z\right)  \label{eq-nablaXeta-star-2}
\end{equation}
Using (\ref{eq-nablaXeta-star-2}) in (\ref{eq-nablaXeta-star-1}), we get (\ref{eq-nablaXeta-star}). 
$\blacksquare$

\begin{thm}
{\rm \cite[Kazan 2019, Theorem~2]{Kazan-19-PhyA}} Let $({\sc D},g)$ be a statistical structure and let $(\varphi,\xi,\eta,g)$ be an almost contact structure on a smooth manifold $M$. Then $(\varphi,\xi,\eta,g,{\sc D})$ is a $(\alpha,\beta)$ trans-Sasakian statistical structure if and only if
\begin{equation}
{\sc D}_{X}\varphi Y-\varphi {\sc D}_{X}^{\ast}Y=\alpha \left(g\left(X,Y\right) \xi -\eta \left(Y\right) X\right) +\beta \left(g\left(\varphi {X},{Y}\right) \xi -\eta \left(Y\right) \varphi {X}\right), 
\label{eq-ts-sta-mfd-1}
\end{equation}
\begin{equation}
{\sc D}_{X}\xi =-\alpha \varphi X+\beta \left(X-\eta (X)\xi \right) +\eta
\left({\sc D}_{X}\xi \right) \xi   \label{eq-ts-sta-mfd-2}
\end{equation}
for all vector fields $X$, $Y$ on $M$.
\end{thm}

\begin{prop}
Let $(M,\varphi,\xi,\eta,g,{\sc D})$ be a trans-Sasakian statistical manifold of type $(\alpha,\beta)$. Then
\begin{equation}
{\sc D}_{X}^{\ast}\xi = -\alpha \varphi {X} + \beta \left(X-\eta \left(X\right) \xi \right) -\eta \left({\sc D}_{X}\xi \right) \xi,
\label{eq-nablaXxi-star-tssm}
\end{equation}
\begin{equation}
\left({\sc D}_{X}\eta \right) Y=-\alpha g(\varphi X,Y)+\beta \left(g(X,Y)-\eta (X)\eta (Y)\right) -\eta \left({\sc D}_{X}\xi \right) \eta (Y),
\label{eq-nablaXeta-tsm-stat}
\end{equation}
\begin{equation}
\left({\sc D}_{X}^{\ast}\eta \right) Y=-\alpha g\left(Y,\varphi X\right)
+\beta \left(g(X,Y)-\eta (X)\eta (Y)\right) +\eta \left({\sc D}_{X}\xi\right) \eta (Y)  \label{eq-nablaXeta-tsm-stat-star}
\end{equation}
for all vector fields $X$, $Y$ on $M$.
\end{prop}

\noindent {\bf Proof.} From (\ref{eq-nabla-star-zero}), we have
\begin{equation}
{\sc D}_{X}\xi +{\sc D}_{X}^{\ast}\xi =2\nabla_{X}\xi.\label{eq-nablaXeta-3}
\end{equation}
Using (\ref{eq-nabla_X(xi)}) and (\ref{eq-ts-sta-mfd-2}) in (\ref{eq-nablaXeta-3}), we obtain (\ref{eq-nablaXxi-star-tssm}). From (\ref{eq-nablaXeta-2}) and (\ref{eq-nablaXxi-star-tssm}), we get (\ref{eq-nablaXeta-tsm-stat}). Also, from (\ref{eq-nablaXeta-star}) and (\ref{eq-ts-sta-mfd-2}), we obtain (\ref{eq-nablaXeta-tsm-stat-star}).
$\blacksquare$

Consider a $(2n+1)$-dimensional almost contact metric statistical manifold $(M,\varphi,\xi,\eta,g,{\sc D})$. A $2$-plane section in $T_{p}M$ is known to be a $\varphi$-section if there is a vector $X\in T_{p}M$ orthogonal to $\xi$ such that $\{X,\varphi {X}\}$ spans the section. The statistical sectional curvature $K^{\sc s}(X,\varphi {X})$, denoted by $c^{\sc s}(X)$, is called a statistical $\varphi$-sectional curvature, given by
\begin{equation}
c^{\sc s}(X) = K^{\sc s}(X,\varphi {X}) = \frac{R^{\sc s}(X,\varphi {X},\varphi{X},X)}{g(X,X)g(\varphi{X},\varphi 
{X})-g(X,\varphi {X})^{2}} = \frac{R^{\sc s}(X,\varphi {X},\varphi {X},X)}{g(X,X)^{2}}.  
\label{eq-stat-phi-sec-cur}
\end{equation}

\begin{defn}
Let $(M,\varphi,\xi,\eta,g,{\sc D})$ be a trans-Sasakian statistical manifold of type $(\alpha,\beta)$, and $R^{\sc s}$ be the statistical curvature tensor on $M$. Then $(M,\varphi,\xi,\eta,g,{\sc D})$ is said to be of pointwise constant statistical $\varphi$-sectional curvature $c^{\sc s}$ if and only if
\begin{eqnarray}
R^{\sc s} &=&\frac{c^{\sc s}+3\left({\alpha}^{2}-{\beta}^{2}\right)}{8}g\circledast g+\frac{c^{\sc s}-\left({\alpha}^{2}-{\beta}^{2}\right)}{4}\left(\frac{1}{2}{\varphi}^{\flat}
\circledast {\varphi}^{\flat} -{\varphi}^{\flat} \circledcirc {\varphi}^{\flat} \right)   \nonumber \\
&&+\left({\rm d}\beta \left(\xi \right) -\frac{c^{\sc s}-\left({\alpha}^{2}-\beta^{2}\right)}{4}\right) g\circledast \left(\eta \otimes \eta\right)   \nonumber \\
&&+g\circledast \left(\left(\varphi {\rm grad}\left(\alpha \right)
\right)^{\sharp}\otimes \eta +\eta \otimes \left(\varphi {\rm grad}\left(\alpha \right) \right)^{\sharp}-\left({\rm d}\beta \otimes \eta +\eta\otimes d\beta \right) \right).  \label{eq-sta-TSSF}
\end{eqnarray}
\end{defn}

\begin{defn-new}
A trans-Sasakian statistical manifold of type $(\alpha,\beta)$ with pointwise constant statistical $\varphi$-sectional curvature $c^{\sc s}$ is called trans-Sasakian statistical space form.
\end{defn-new}

\begin{rem-new}
By substituting the particular values of $\alpha$ and $\beta$ corresponding to the subclasses listed in Remark \ref{rem-class-tssm}, into $(\ref{eq-sta-TSSF})$, one obtains the explicit expressions for Riemann--Christoffel curvature tensors of ${\cal C}_{6}$-statistical space forms, $\alpha$-Sasakian statistical space forms, Sasakian statistical space forms, ${\cal C}_{5}$-statistical space forms, $\beta$-Kenmotsu statistical space forms, Kenmotsu statistical space forms and cosymplectic statistical space forms.
\end{rem-new}

\begin{ex-new}
\label{ex-S=R0}Let $(M,\varphi,\xi,\eta,g,\nabla)$ be a $(2n+1)$-dimensional almost contact metric manifold. Define difference tensor field ${\sc B}^{\sc d}$ by
\begin{equation}
{\sc B}^{\sc d}_{X}Y = f\eta \left(X\right) \eta \left(Y\right) \xi 
\label{eq-trivial-statistical-structure}
\end{equation}
for all vector fields $X$, $Y$ on $M$, where $f$ is a smooth function on $M$. It is straightforward to verify that ${\sc B}^{\sc d}$ satisfies conditions (\ref{eq-K-properties}) and (\ref{eq-almost-con-statis}). Consequently $(\varphi,\xi,\eta,g,{\sc D}=\nabla+{\sc B}^{\sc d})$ is an almost contact metric statistical structure on $M$. Moreover, one easily checks that ${\sc B}^{\sc d}$ also satisfies
\[
[{\sc B}^{\sc d}_{X},{\sc B}^{\sc d}_{Y}]Z=0
\]
for all vector fields $X$, $Y$, $Z$ on $M$. Thus, from (\ref{eq-statistical-cur}), we have $R^{\sc s}=R$. In particular, if $M$ is a trans-Sasakian space form whose $\varphi$-sectional curvature is pointwise constant $c$, then $(M,\varphi,\xi,\eta,g,{\sc D} =\nabla+{\sc B}^{\sc d})$ is a trans-Sasakian statistical space form with pointwise constant statistical $\varphi$-sectional curvature $c^{\sc s}=c$.
\end{ex-new}

\begin{rem-new}
The structure given by (\ref{eq-trivial-statistical-structure}) on an almost contact metric statistical manifold is called almost-trivial for the following two reasons \cite[Wu {\em et al.} 2024, Remark~2.4]{Wu-Sheng-Zhang-24-DGA}. First, its statistical curvature tensor coincides with the Riemannian curvature tensor. Second, the difference tensor ${\sc B}^{\sc d}$ satisfies
\[
{\sc B}^{\sc d}(X,Y) = 0
\]
whenever either $X$ or $Y$ is orthogonal to $\xi$.
\end{rem-new}

\section{Trans-Sasakian statistical space forms of dimension  
$>3$\label{sec-tran-Sasa-Stat-nonproper}}

The purpose of this section is to prove that every trans-Sasakian statistical space form of dimension greater than $3$ is necessarily almost-trivial.

\begin{lem}
\label{Lemma-almostContactStatistical} {\rm \cite[Wu {\em et al.} 2024]{Wu-Sheng-Zhang-24-DGA}} Let $(M,\varphi,\xi,\eta,g,{\sc D} =\nabla+{\sc B}^{\sc d})$ be an almost contact metric statistical manifold. Then,
\begin{equation}
{\sc B}^{\sc d}_{X}\xi = \eta\left({\sc B}^{\sc d}_{\xi}\xi\right)\eta\left(X\right)\xi,  \label{eq-K(X,xi)-ACMS}
\end{equation}
\begin{equation}
\eta \left({\sc B}^{\sc d}_{X}Y \right) =\eta \left({\sc B}^{\sc d}_{\xi}\xi \right) \eta \left(X\right) \eta \left(Y\right)  \label{eq-etaK(X,Y)-ACMS}
\end{equation}
for all vector fields $X$, $Y$ on $M$.
\end{lem}

Let $(M,\nabla,\varphi,\xi,\eta,g)$ be a $(2n+1)$-dimensional trans-Sasakian space form with pointwise constant $\varphi$-sectional curvature $c$, and let $(M,\varphi,\xi,\eta,g,{\sc D} =\nabla+{\sc B}^{\sc d})$ be a trans-Sasakian statistical space form with pointwise constant statistical $\varphi$-sectional curvature $c^{\sc s}$. It follows from (\ref{eq-statistical-cur}), (\ref{eq-tps-Riem-cur-for-33}) and (\ref{eq-sta-TSSF}) that
\begin{equation}
g\left([{\sc B}^{\sc d}_{X},{\sc B}^{\sc d}_{Y}]Z,W\right) =\frac{c^{\sc s}-c}{8}\{g\circledast g+{\varphi}^{\flat} \circledast {\varphi}^{\flat} -2{\varphi}^{\flat} \circledcirc {\varphi}^{\flat} -2g\circledast \left(\eta \otimes \eta \right) \}\left(X,Y,Z,W\right) .  \label{eq-KxKy-TSSF}
\end{equation}
Let $\{e_{1},\ldots,e_{n},\varphi e_{1},\ldots,\varphi e_{n},\xi \}$ be a $\varphi $-basis of $M$. Then from the Lemma \ref{Lemma-almostContactStatistical}, we get
\begin{equation}
{\sc B}^{\sc d}_{e_{i}}\xi=0,\quad \eta ({\sc B}^{\sc d}_{e_{i}}e_{j})=0
\label{eq-K(ei,xi)-etaK(ei,ej)-1}
\end{equation}
for $i,j\in \{1,\ldots ,n\}$. From (\ref{eq-K-properties}), (\ref{eq-almost-con-statis}), (\ref{eq-K(ei,xi)-etaK(ei,ej)-1}), we obtain
\begin{eqnarray*}
{\sc B}^{\sc d}_{e_{i}}\varphi e_{j} &=& {\sc B}^{\sc d}_{\varphi e_{i}}e_{j} = -\varphi {\sc B}^{\sc d}_{e_{i}}e_{j},\\
{\sc B}^{\sc d}_{\varphi e_{i}}\varphi e_{j} &=&\varphi^{2}{\sc B}^{\sc d}_{e_{i}}e_{j} = -{\sc B}^{\sc d}_{e_{i}}e_{j}.
\end{eqnarray*}
Also, we have
\begin{eqnarray*}
g(\varphi{\sc B}^{\sc d}_{e_{i}}e_{j},\varphi {\sc B}^{\sc d}_{e_{k}}e_{\ell}) &=& g({\sc B}^{\sc d}_{e_{i}}e_{j},{\sc B}^{\sc d}_{e_{k}}e_{\ell}),\\
g(\varphi {\sc B}^{\sc d}_{e_{i}}e_{j},{\sc B}^{\sc d}_{e_{k}}e_{\ell})
 &=& -g\left({\sc B}^{\sc d}_{e_{i}}e_{j}, \varphi {\sc B}^{\sc d}_{e_{k}}e_{\ell}\right),
\end{eqnarray*}
where $i,j,k,\ell \in \{1,\ldots ,n\}$.

\begin{lem}
\label{Lemma-K(ei,ej)-cond}Let $(M,\varphi,\xi,\eta,g,\nabla)$ be a $(2n+1)$-dimensional trans-Sasakian space form with $n>1$ and pointwise constant $\varphi$-sectional curvature $c$. Assume that $(\varphi,\xi,\eta,g,{\sc D}=\nabla+{\sc B}^{\sc d})$ admits a trans-Sasakian statistical structure on $M$. If $(M,\varphi,\xi,\eta,g,{\sc D})$ is a trans-Sasakian statistical space form with pointwise constant statistical $\varphi$-sectional curvature $c^{\sc s}$. Then for $i\neq j$, 
\begin{equation}
g ({\sc B}^{\sc d}_{e_{i}}e_{i},{\sc B}^{\sc d}_{e_{i}}e_{i}) = \frac{c-c^{\sc s}}{2},  \label{eq-K(ei,ej)-1}
\end{equation}
\begin{equation}
g({\sc B}^{\sc d}_{e_{i}}e_{j},{\sc B}^{\sc d}_{e_{i}}e_{j}) = \frac{c-c^{\sc s}}{4},  \label{eq-K(ei,ej)-2}
\end{equation}
\begin{equation}
g\left({\sc B}^{\sc d}_{e_{j}}e_{j},{\sc B}^{\sc d}_{e_{i}}e_{i}\right) =0,
\label{eq-K(ei,ej)-3}
\end{equation}
\begin{equation}
g\left({\sc B}^{\sc d}_{e_{j}}e_{j},\varphi {\sc B}^{\sc d}_{e_{i}}e_{i}\right) =0.  \label{eq-K(ei,ej)-4}
\end{equation}
Also, for $j\neq i,k$ 
\begin{equation}
g\left({\sc B}^{\sc d}_{e_{k}}e_{j},{\sc B}^{\sc d}_{e_{i}}e_{i}\right) =0,
\label{eq-K(ei,ej)-5}
\end{equation}
\begin{equation}
g\left({\sc B}^{\sc d}_{e_{k}}e_{j},\varphi{\sc B}^{\sc d}_{e_{i}}e_{i}\right) = 0.  \label{eq-K(ei,ej)-6}
\end{equation}

\end{lem}

\noindent {\bf Proof.} Putting $X=Z=e_{i}$ and $Y=W=\varphi e_{i}$ in (\ref{eq-KxKy-TSSF}), and using (\ref{eq-almost-con-statis}) and (\ref{eq-K(ei,xi)-etaK(ei,ej)-1}) in the resulting equation, we get (\ref{eq-K(ei,ej)-1}). Putting $X=Z=e_{i}$ and $Y=W=e_{j}$, $i\neq j$, in (\ref{eq-KxKy-TSSF}), we get 
\begin{equation}
g\left({\sc B}^{\sc d}_{e_{i}}e_{j},{\sc B}^{\sc d}_{e_{i}}e_{j}\right) - g\left({\sc B}^{\sc d}_{e_{j}}e_{j},{\sc B}^{\sc d}_{e_{i}}e_{i}\right) =\frac{c-c^{\sc s}}{4}.  \label{eq-gK(ei,ej)-gK(ej,ej)}
\end{equation}
Replacing $e_{i}$ by $\varphi e_{i}$ in (\ref{eq-gK(ei,ej)-gK(ej,ej)}), we obtain
\begin{equation}
g\left({\sc B}^{\sc d}_{e_{i}}e_{j},{\sc B}^{\sc d}_{e_{i}}e_{j}\right)
+g\left({\sc B}^{\sc d}_{e_{j}}e_{j},{\sc B}^{\sc d}_{e_{i}}e_{i} \right) = \frac{c-c^{\sc s}}{4}.  \label{eq-gK(ei,ej)-gK(ej,ej)-1}
\end{equation}
From (\ref{eq-gK(ei,ej)-gK(ej,ej)}) and (\ref{eq-gK(ei,ej)-gK(ej,ej)-1}), we get
\[
g\left({\sc B}^{\sc d}_{e_{i}}e_{j},{\sc B}^{\sc d}_{e_{i}}e_{j}\right) =\frac{c-c^{\sc s}}{4}, 
\]
\[
g\left({\sc B}^{\sc d}_{e_{j}}e_{j},{\sc B}^{\sc d}_{e_{i}}e_{i}\right) =0, 
\]
which gives (\ref{eq-K(ei,ej)-2}) and (\ref{eq-K(ei,ej)-3}).
Putting $X=e_{i}$, $Y=W=e_{j}$, $Z=\varphi e_{i}$ in (\ref{eq-KxKy-TSSF}),
we obtain
\begin{equation}
g\left({\sc B}^{\sc d}_{e_{i}}e_{j},\varphi{\sc B}^{\sc d}_{e_{i}}e_{j}\right) - g\left({\sc B}^{\sc d}_{e_{j}}e_{j},\varphi {\sc B}^{\sc d}_{e_{i}}e_{i}\right) =0.  \label{eq-gK(ei,ej)-gK(ej,ej)-2}
\end{equation}
Since $g\left({\sc B}^{\sc d}_{e_{i}}e_{j},\varphi {\sc B}^{\sc d}_{e_{i}}e_{j}\right) =0$, therefore (\ref{eq-K(ei,ej)-4}) follows from (\ref{eq-gK(ei,ej)-gK(ej,ej)-2}).
Putting $X=e_{k}$, $Y=Z=\varphi e_{i}$, $W=e_{j}$ in (\ref{eq-KxKy-TSSF}), we get
\begin{equation}
g\left({\sc B}^{\sc d}_{e_{k}}e_{j},{\sc B}^{\sc d}_{e_{i}}e_{i}\right) + g\left({\sc B}^{\sc d}_{e_{i}}e_{j},{\sc B}^{\sc d}_{e_{k}}e_{i}\right) = 0.
\label{eq-gK(ei,ej)-gK(ej,ej)-3}
\end{equation}
Replacing $e_{i}$ by $\varphi e_{i}$ in (\ref{eq-gK(ei,ej)-gK(ej,ej)-3}), we get
\begin{equation}
g\left({\sc B}^{\sc d}_{e_{k}}e_{j},{\sc B}^{\sc d}_{e_{i}}e_{i}\right) = g\left({\sc B}^{\sc d}_{e_{i}}e_{j},{\sc B}^{\sc d}_{e_{k}}e_{i}\right).
\label{eq-gK(ei,ej)-gK(ej,ej)-4}
\end{equation}
Hence from (\ref{eq-gK(ei,ej)-gK(ej,ej)-3}) and (\ref{eq-gK(ei,ej)-gK(ej,ej)-4}), we get (\ref{eq-K(ei,ej)-5}). By replacing $e_{j}$ with $\varphi e_{j}$ in (\ref{eq-K(ei,ej)-5}), we obtain (\ref{eq-K(ei,ej)-6}). $\blacksquare$

\begin{thm}\label{th-tsm-dim>3}
Let $(M,\varphi,\xi,\eta,g,\nabla)$ be a $(2n+1)$-dimensional trans-Sasakian space form with $n>1$ and pointwise constant $\varphi$-sectional curvature $c$. Assume that $(M,\varphi,\xi,\eta,g,{\sc D} =\nabla+{\sc B}^{\sc d})$ admits a trans-Sasakian statistical structure on $M$. If $(M,\varphi,\xi,\eta,g,{\sc D})$ is a trans-Sasakian statistical space form with pointwise constant statistical $\varphi$-sectional curvature $c^{\sc s}$. Then,
\[
c^{\sc s}=c,
\]
and the difference tensor field ${\sc B}^{\sc d}$ is given by
\[
{\sc B}^{\sc d}_{X}Y = f\eta \left(X\right) \eta \left(Y\right) \xi 
\]
for all vector fields $X$, $Y$ on $M$, where $f$ is a smooth function on $M$.
\end{thm}

\noindent {\bf Proof.} Assume that $c\neq c^{\sc s}$. Then, from (\ref{eq-K(ei,ej)-2}), it follows that
\begin{equation}
{\sc B}^{\sc d}_{e_{1}}e_{2} \neq 0.  \label{eq-K(e1,e2)neq0}
\end{equation}
Let $\{e_{1},\ldots ,e_{n},\varphi e_{1},\ldots,\varphi e_{n},\xi \}$ be a $\varphi$-basis of $M$. By virtue of (\ref{eq-K(ei,ej)-1}) and (\ref{eq-K(ei,ej)-3}), the vectors ${\sc B}^{\sc d}_{e_{i}}e_{i}$, $i\in\{1,\ldots,n\}$, are non-zero and orthogonal to each other. Also from (\ref{eq-K(ei,xi)-etaK(ei,ej)-1}), (\ref{eq-K(ei,ej)-1}), (\ref{eq-K(ei,ej)-3}), we get
\[
g(\varphi {\sc B}^{\sc d}_{e_{i}}e_{i},\varphi {\sc B}^{\sc d}_{e_{i}}e_{i}) = g({\sc B}^{\sc d}_{e_{i}}e_{i},{\sc B}^{\sc d}_{e_{i}}e_{i})=\frac{c-c^{\sc s}}{2}, 
\]
and
\[
g\left(\varphi {\sc B}^{\sc d}_{e_{i}}e_{i}, \varphi {\sc B}^{\sc d}_{e_{j}}e_{j}\right) = g\left({\sc B}^{\sc d}_{e_{i}}e_{i},{\sc B}^{\sc d}_{e_{j}}e_{j} \right) =0,\quad i\neq j. 
\]
From (\ref{eq-K(ei,xi)-etaK(ei,ej)-1}), it is clear that ${\sc B}^{\sc d}_{e_{i}}e_{i}$ is orthogonal to $\xi$. Since 
\[
g\left(\varphi {\sc B}^{\sc d}_{e_{i}}e_{i},\xi \right) = - g\left({\sc B}^{\sc d}_{e_{i}}e_{i},\varphi\xi\right) =0,\quad i\in\{1,\ldots,n\}, 
\]
therefore $\varphi {\sc B}^{\sc d}_{e_{i}}e_{i}$ is also orthogonal to $\xi$. In addition, (\ref{eq-K(ei,ej)-4}) implies that ${\sc B}^{\sc d}_{e_{i}}e_{i}$ is orthogonal to $\varphi {\sc B}^{\sc d}_{e_{j}}e_{j}$ for every $i,j\in \{1,\ldots ,n\}$. Consequently, 
\[
\{{\sc B}^{\sc d}_{e_{i}}e_{i},\varphi {\sc B}^{\sc d}_{e_{i}}e_{i},\xi \},\quad i\in \{1,\ldots ,n\}, 
\]
forms an orthogonal basis of $M$. It follows from Lemma \ref{Lemma-K(ei,ej)-cond} that
\[
g({\sc B}^{\sc d}_{e_{1}}e_{2},{\sc B}^{\sc d}_{e_{i}}e_{i}) = 0,\quad g({\sc B}^{\sc d}_{e_{1}}e_{2},\varphi {\sc B}^{\sc d}_{e_{i}}e_{i})=0,\quad i\in \{1,\ldots ,n\}. 
\]
Moreover, by (\ref{eq-K(ei,xi)-etaK(ei,ej)-1}), we have 
\[
g({\sc B}^{\sc d}_{e_{1}}e_{2},\xi)=0. 
\]
Hence ${\sc B}^{\sc d}_{e_{1}}e_{2}=0$, which contradicts (\ref{eq-K(e1,e2)neq0}). Therefore, we conclude that $c=c^{\sc s}$. We now describe the associated statistical structure. Since $c=c^{\sc s}$, it follows from (\ref{eq-g(phiX,phiY)}), (\ref{eq-almost-con-statis}), (\ref{eq-etaK(X,Y)-ACMS}), (\ref{eq-K(ei,ej)-1}) and (\ref{eq-K(ei,ej)-2}), that
\begin{equation}
{\sc B}^{\sc d}_{e_{i}}e_{j} = {\sc B}^{\sc d}_{e_{i}}\varphi e_{j} = {\sc B}^{\sc d}_{\varphi e_{i}}\varphi e_{j}=0,\quad i,j\in \{1,\ldots ,n\}.
\label{eq-K(ei,ej)=K(ei,phiej)}
\end{equation}
Replacing $X$ by $\xi$, $e_{i}$ and $\varphi e_{i}$, respectively, in (\ref{eq-K(X,xi)-ACMS}), we obtain
\begin{equation}
{\sc B}^{\sc d}_{\xi}\xi = \eta \left({\sc B}^{\sc d}_{\xi}\xi\right) \xi,\quad {\sc B}^{\sc d}_{e_{i}}\xi = {\sc B}^{\sc d}_{\varphi e_{i}}\xi =0,\quad i\in \{1,\ldots ,n\}.  \label{eq-K(xi,xi)=eta(K(xi,xi)xi)}
\end{equation}
Suppose
\begin{equation}
X=\sum\limits_{i=1}^{n}(\lambda_{i}e_{i}+\lambda_{n+i}\varphi e_{i})+\lambda_{2n+1}\xi ,  \label{eq-X}
\end{equation}
and
\begin{equation}
Y=\sum\limits_{j=1}^{n}(\mu_{j}e_{j}+\mu_{n+j}\varphi e_{j})+\mu_{2n+1}\xi .  \label{eq-Y}
\end{equation}
Then, by using (\ref{eq-K(ei,ej)=K(ei,phiej)}) and (\ref{eq-K(xi,xi)=eta(K(xi,xi)xi)}), we find
\[
{\sc B}^{\sc d}_{X}Y = \lambda_{2n+1}\mu_{2n+1}\eta \left({\sc B}^{\sc d}_{\xi}\xi \right)\xi.
\]
Finally, from (\ref{eq-X}) and (\ref{eq-Y}), we have $\lambda_{2n+1}=\eta(X)$, $\mu_{2n+1}=\eta (Y)$. Setting $f=\eta ({\sc B}^{\sc d}_{\xi}\xi)$, we conclude that ${\sc B}^{\sc d}(X,Y)=f\eta (X)\eta (Y)\xi $. $\blacksquare$

\begin{rem-new}
One can readily verify that Theorem~\ref{th-tsm-dim>3} remains valid for the subclasses of trans-Sasakian statistical manifolds described in Remark~\ref{rem-class-tssm}.
\end{rem-new}

\begin{ex-new}
Let $M={\Bbb R}^{5}$ be a $5$-dimensional smooth manifold with standard coordinates $\left(x^{1},x^{2},x^{3},x^{4},x^{5}\right)$. Define a linearly independent global coordinate frame $\{U_{1},U_{2},U_{3},U_{4},U_{5}\}$ on $M$ by
\[
U_{i}=e^{-x^{5}}\frac{\partial }{\partial x^{i}},\quad i\in \{1,\ldots ,5\}. 
\]
Consider the Riemannian metric $g$ on $M$ defined by
\[
g(U_{i},U_{j})= \delta_{ij}, \;\; i,j\in \{1,\dots,5\},
\]
where $\delta_{ij}$ is the well-known Kronecker delta.
Let $U_{5}=\xi$, and define the $1$-form $\eta$ on $M$ by
\[
\eta (X)=g(X,U_{5}) 
\]
for every vector field $X$ on $M$. Also, Let $\varphi$ be the $(1,1)$ tensor field given by
\[
\varphi (U_{1})=U_{3},\quad \varphi (U_{2})=U_{4},\quad \varphi
(U_{3})=-U_{1},\quad \varphi (U_{4})=-U_{2},\quad \varphi (U_{5})=0. 
\]
Let $\nabla$ denote the Levi--Civita connection of $g$, and let $R$ be the Riemann--Christoffel curvature tensor of $M$. Then, by a direct calculation, we obtain 
\begin{eqnarray*}
\lbrack U_{i},U_{j}] &=&0,\quad i,j\in \{1,\ldots ,4\},\\
\lbrack U_{i},U_{5}] &=&e^{-x^{5}}U_{i,}\quad i\in \{1,\ldots ,4\}.\\
\end{eqnarray*}
Also, we have 
\begin{eqnarray*}
{\nabla}_{U_{i}}U_{i} &=&-e^{-x^{5}}U_{i},\quad i\in \{1,\ldots ,4\},\\
{\nabla}_{U_{i}}U_{j} &=&0,\quad 1\leq i\neq j\leq 4,\\
{\nabla}_{U_{5}}U_{i}&=&0, \quad{\nabla}_{U_{i}}U_{5}=e^{-x^{5}}U_{i},\;\; i\in\{1,\ldots,4\},\\
{\nabla}_{U_{5}}U_{5} &=&0.
\end{eqnarray*}
Then, $(\varphi,\xi,\eta,g)$ defines a trans Sasakian structure on $M$ with $\left(\alpha,\beta\right) =\left(0,e^{-x^{5}}\right) $ (see \cite[Hu and Wang 2016]{Hu-Wang-16-IEJG}).
The required components of the Riemann--Christoffel curvature tensor $R$ are obtained from the above relations as follows:
\begin{eqnarray*}
R(U_{1},U_{3})U_{3} &=& -e^{-2x^{5}}U_{1}, \qquad R(U_{2},U_{4})U_{4} = -e^{-2x^{5}}U_{2}\\
R(U_{3},U_{1})U_{1} &=& -e^{-2x^{5}}U_{3}, \qquad R(U_{4},U_{2})U_{2} = -e^{-2x^{5}}U_{4}.
\end{eqnarray*}
The $\varphi$-sectional curvature of $M$, computed from the curvature tensor components, is
\[
K(U_{1},\varphi U_{1}) = K(U_{2},\varphi
U_{2}) = K(U_{3},\varphi U_{3}) = K(U_{4},\varphi
U_{4})=-e^{-2x^{5}}. 
\]
Hence $M$ is a trans-Sasakian space form with pointwise constant $\varphi$-sectional curvature $-e^{-2x^{5}}$. Define an affine connection ${\sc D}$ on $M$ as
\begin{eqnarray*}
{\sc D}_{U_{i}}U_{i} &=&-e^{-x^{5}}U_{i},\quad i\in \{1,\ldots,4\},\\
{\sc D}_{U_{i}}U_{j} &=&0,\quad 1\leq i\neq j\leq 4,\\
{\sc D}_{U_{5}}U_{i}&=&0,\quad
{\sc D}_{U_{i}}U_{5}=e^{-x^{5}}U_{i},\quad i\in\{1,\ldots,4\},\\
{\sc D}_{U_{5}}U_{5} &=&\lambda \xi ,
\end{eqnarray*}
where $\lambda $ is a real number. Since $K={\sc D} -\nabla^{0}$, then we obtain
\begin{eqnarray*}
{\sc B}^{\sc d}_{U_{i}}U_{i} &=& 0,\quad i\in \{1,\ldots ,4\},\\
{\sc B}^{\sc d}_{U_{i}}U_{j} &=&0,\quad 1\leq i\neq j\leq5,\\
{\sc B}^{\sc d}_{U_{5}}U_{5}&=&\lambda \xi
\end{eqnarray*}
or
\[
{\sc B}^{\sc d}_{X}Y = \lambda \eta \left(X\right) \eta \left(Y\right) \xi 
\]
for all vector fields $X$, $Y$ on $M$. Hence in view of Example~\ref{ex-S=R0}, $(M,\varphi,\xi,\eta,g,{\sc D})$ is a trans-Sasakian statistical space form with pointwise constant statistical $\varphi$-sectional curvature $-e^{-2x^{5}}$.
\end{ex-new}

\section{$3$-dimensional trans-Sasakian statistical space forms\label{sec-3D-trans-sasa-stat}}

In this section, we consider trans-Sasakian statistical structures on $3$-diemnsional trans-Sasakian space forms. We prove that such manifolds admit non-trivial trans-Sasakian statistical structures with pointwise constant statistical $\varphi$-sectional curvature. This contrasts with the higher-dimensional case, where every such structure is necessarily almost-trivial. 

\begin{thm}\label{th-tsm-dim=3}
\label{th-tssf-sec-cur-inequality} Let $(M,\varphi,\xi,\eta,g,\nabla)$ be a $3$-dimensional trans-Sasakian space form with pointwise constant $\varphi$-sectional curvature $c$, and let $(M,\varphi,\xi,\eta,g,{\sc D} =\nabla+{\sc B}^{\sc d})$ be a trans-Sasakian statistical space form with pointwise constant statistical $\varphi$-sectional curvature $c^{\sc s}$. Then,
\[
c^{\sc s}\leq c. 
\]
Further, with respect to an orthonormal frame $\{e_{1},\varphi e_{1},\xi \}$, the difference tensor ${\sc B}^{\sc d}$ is given by
\begin{eqnarray*}
{\sc B}^{\sc d}_{e_{1}}e_{1} &=& f_{1}e_{1}+f_{2}\varphi e_{1},\quad {\sc B}^{\sc d}_{\varphi e_{1}}\varphi e_{1}=-f_{1}e_{1}-f_{2}\varphi e_{1},\\
{\sc B}^{\sc d}_{e_{1}}\varphi e_{1} &=&{\sc B}^{\sc d}_{\varphi e_{1}}e_{1} = f_{2}e_{1}-f_{1}\varphi e_{1},\quad {\sc B}^{\sc d}_{\xi}\xi =f\xi ,\\
{\sc B}^{\sc d}_{e_{1}}\xi &=& {\sc B}^{\sc d}_{\xi} e_{1} = {\sc B}^{\sc d}_{\varphi e_{1}}\xi  = {\sc B}^{\sc d}_{\xi},\varphi e_{1} =0,
\end{eqnarray*}
where $f_{1}$, $f_{2}$, $f$ are smooth functions on $M$ and satisfy
\[
f_{1}^{2}+f_{2}^{2}=\frac{c-c^{\sc s}}{2}. 
\]
\end{thm}

\noindent {\bf Proof.} Let $\{e_{1},e_{2}=\varphi e_{1},e_{3}=\xi\}$ be an orthonormal basis of $(M,\varphi,\xi,\eta,g,{\sc D} =\nabla+{\sc B}^{\sc d})$. Then,
\begin{equation}
{\sc B}^{\sc d}_{e_{i}}e_{j} = {\sc B}^{\sc d}_{e_{j}}e_{i},\quad g\left({\sc B}^{\sc d}_{e_{i}}e_{j},e_{k}\right) = g({\sc B}^{\sc d}_{e_{i}}e_{k},e_{j}), \quad {\sc B}^{\sc d}_{e_{i}}\varphi e_{j}+\varphi {\sc B}^{\sc d}_{e_{i}}e_{j} =0.
\label{eq-K(ei,ej)=K(ej,ei)}
\end{equation}
where $i,j,k\in \{1,2,3\}$. By virtue of (\ref{eq-K(ei,xi)-etaK(ei,ej)-1}), we may suppose that
\begin{equation}
{\sc B}^{\sc d}_{e_{1}}e_{1} = f_{1}e_{1}+f_{2}e_{2}.\label{eq-K(e1,e1)=lambda1e1+lambda2e2}
\end{equation}
From (\ref{eq-K(X,xi)-ACMS}) and (\ref{eq-K(ei,ej)=K(ej,ei)}), we have
\begin{equation}
\left\{
\begin{array}{l}
{\sc B}^{\sc d}_{e_{2}}e_{2}=-f_{1}e_{1}-f_{2}e_{2},\quad
{\sc B}^{\sc d}_{e_{1}}e_{2}={\sc B}^{\sc d}_{e_{2}}e_{1}=-f_{2}e_{1}-f_{1}e_{2}, \\
{\sc B}^{\sc d}_{e_{1}}\xi ={\sc B}^{\sc d}_{\xi}e_{1}={\sc B}^{\sc d}_{e_{2}}\xi={\sc B}^{\sc 
d}_{\xi}e_{2}=0,\quad {\sc B}^{\sc d}_{\xi}\xi=f \xi  
\end{array}
\right. \label{eq-K(e2,e2)=-lambda1e1-lambda2e2} 
\end{equation}
for $f=\eta ({\sc B}^{\sc d}_{\xi}\xi)$.
Noting that (\ref{eq-K(ei,ej)-1}) also holds for dimension $3$, therefore
\[
g({\sc B}^{\sc d}_{e_{1}}e_{1},{\sc B}^{\sc d}_{e_{1}}e_{1}) = \frac{c-c^{\sc s}}{2}, 
\]
which gives $c^{\sc s}\leq c$. By using (\ref{eq-K(e1,e1)=lambda1e1+lambda2e2}) in (\ref{eq-K(e2,e2)=-lambda1e1-lambda2e2}), we get
\[
f_{1}^{2}+f_{2}^{2}=\frac{c-c^{\sc s}}{2}. 
\]
This completes the proof. $\blacksquare$

\begin{rem-new}
One can readily verify that Theorem~\ref{th-tsm-dim=3} remains valid for the subclasses of trans-Sasakian statistical manifolds described in Remark~\ref{rem-class-tssm}.
\end{rem-new}

Next, we illustrate the inequality $c^{\sc s}\leq c$ by means of examples
on $3$-dimensional proper trans-Sasakian space forms.

\begin{ex-new}
Let $M=\{\left(x^{1},x^{2},x^{3}\right) \in {\Bbb R}^{3}:x^{3}>0\}$ be a $3$-dimensional smooth manifold with standard coordinates $\left(x^{1},x^{2},x^{3}\right)$. Define a linearly independent global coordinate frame $\{U_{1},U_{2},U_{3}\}$ on $M$ by
\[
U_{1}=x^{3}\left(\frac{\partial}{\partial x^{1}}+x^{2}\frac{\partial}{\partial x^{3}}\right) ,\quad U_{2}=x^{3}\frac{\partial }{\partial x^{2}},\quad U_{3}=\frac{\partial}{\partial x^{3}}. 
\]
Consider the Riemannian metric $g$ on $M$ defined by
\[
g(U_{i},U_{j})= \delta_{ij}, \;\; i,j\in \{1,2,3\},
\]
where $\delta_{ij}$ is the well-known Kronecker delta. 
Let $U_{3}=\xi$, and define the $1$-form $\eta $ on $M$ by
\[
\eta (X)=g(X,U_{3})
\]
for every vector field $X$ on $M$. Also, let $\varphi$ be the $(1,1)$
tensor field given by
\[
\varphi (U_{1})=U_{2},\quad \varphi (U_{2})=-U_{1},\quad \varphi (U_{3})=0. 
\]
Let $\nabla$ denote the Levi--Civita connection of $g$, and let $R$ be the Riemann--Christoffel curvature tensor of $M$. Then, by a direct calculation, we obtain 
\[
[U_{1},U_{2}]=x^{2}U_{2}-\left(x^{3}\right)^{2}U_{3},\quad [U_{1},U_{3}]=-\frac{1}{x^{3}}U_{1},\quad [U_{2},U_{3}]=-\frac{1}{x^{3}}U_{2}. 
\]
Also, we obtain 
\begin{equation*}
\begin{aligned}
\nabla_{U_{1}}U_{1} &= \frac{1}{x^{3}}U_{3}, 
&\qquad \nabla_{U_{1}}U_{2} &= -\frac{1}{2}(x^{3})^{2}U_{3},
&\qquad \nabla_{U_{1}}U_{3} &= -\frac{1}{x^{3}}U_{1}+\frac{1}{2}(x^{3})^{2}U_{2}, \\[0.6em]
\nabla_{U_{2}}U_{1} &= \frac{1}{2}(x^{3})^{2}U_{3}-x^{2}U_{2},
&\qquad \nabla_{U_{2}}U_{2} &= x^{2}U_{1}+\frac{1}{x^{3}}U_{3},
&\qquad \nabla_{U_{2}}U_{3} &= -\frac{1}{x^{3}}U_{2}-\frac{1}{2}(x^{3})^{2}U_{1}, \\[0.6em]
\nabla_{U_{3}}U_{1} &= \frac{1}{2}(x^{3})^{2}U_{2},
&\qquad \nabla_{U_{3}}U_{2} &= -\frac{1}{2}(x^{3})^{2}U_{1},
&\qquad \nabla_{U_{3}}U_{3} &= 0.
\end{aligned}
\end{equation*}
Then, $(\varphi,\xi,\eta,g)$ defines a trans Sasakian structure on $M$ with $\left(\alpha,\beta\right) =\left(-\frac{1}{2}{(x^{3})}^{2},-\frac{1}{x^{3}}\right)$ (see \cite[Shaikh {\em et al.} 2008, Example~2]{Shaikh-Baishya-Eyasmin-08-DM}). The components of the Riemann--Christoffel curvature tensor $R$ are obtained from the above relations as follows.
\begin{eqnarray*}
R(U_{1},U_{2})U_{1}
&=&\left(\frac{3}{4}(x^{3})^{4}
+\frac{1}{(x^{3})^{2}}
+(x^{2})^{2}\right)U_{2}
+x^{2}(x^{3})^{2}U_{3},\\
R(U_{1},U_{2})U_{2}
&=&-\left(\frac{3}{4}(x^{3})^{4}
+\frac{1}{(x^{3})^{2}}
+(x^{2})^{2}\right)U_{1}
-\frac{x^{2}}{x^{3}}U_{3},\\
R(U_{1},U_{2})U_{3}
&=&-x^{2}(x^{3})^{2}U_{1}
+\frac{x^{2}}{x^{3}}U_{2},\\[0.5em]
R(U_{1},U_{3})U_{1}
&=&x^{2}(x^{3})^{2}U_{2}
+\left(\frac{2}{(x^{3})^{2}}
-\frac{1}{4}(x^{3})^{4}\right)U_{3},\\
R(U_{1},U_{3})U_{2}
&=&-x^{2}(x^{3})^{2}U_{1},\\
R(U_{1},U_{3})U_{3}
&=&\left(\frac{1}{4}(x^{3})^{4}
-\frac{2}{(x^{3})^{2}}\right)U_{1},\\[0.5em]
R(U_{2},U_{3})U_{1}
&=&-\frac{x^{2}}{x^{3}}U_{2},\\
R(U_{2},U_{3})U_{2}
&=&\frac{x^{2}}{x^{3}}U_{1}
+\left(\frac{2}{(x^{3})^{2}}
-\frac{1}{4}(x^{3})^{4}\right)U_{3},\\
R(U_{2},U_{3})U_{3}
&=&-\left(\frac{2}{(x^{3})^{2}}
-\frac{1}{4}(x^{3})^{4}\right)U_{2}.
\end{eqnarray*}
The $\varphi$-sectional curvature of $M$, computed from the curvature tensor components, is
\[
c=K\left(U_{1},\varphi U_{1}\right) =K\left(U_{2},\varphi U_{2}\right) =-\left(\frac{3}{4}(x^{3})^{4}+\frac{1}{(x^{3})^{2}}+(x^{2})^{2}\right). 
\]
Hence $M$ is a trans-Sasakian space form with pointwise constant $\varphi$-sectional curvature 
\[-\frac{3}{4}(x^{3})^{4}-\frac{1}{(x^{3})^{2}}-(x^{2})^{2}.
\]
We define an affine connection ${\sc D}$ on $M$ by
\begin{equation*}
\begin{aligned}
{\sc D}_{U_{1}}U_{1} 
&= \frac{1}{x^{3}}U_{3}+\lambda_{1}U_{1}+\lambda_{2}U_{2},
&\qquad {\sc D}_{U_{1}}U_{2} & = -\frac{1}{2}(x^{3})^{2}U_{3}
-\lambda_{1}U_{1}-\lambda_{2}U_{2},\\[0.5em]
{\sc D}_{U_{1}}U_{3} & = -\frac{1}{x^{3}}U_{1}+\frac{1}{2}(x^{3})^{2}U_{2},
&\qquad {\sc D}_{U_{2}}U_{1} & = \frac{1}{2}(x^{3})^{2}U_{3}
-x^{2}U_{2}-\lambda_{1}U_{1}-\lambda_{2}U_{2},\\[0.5em]
{\sc D}_{U_{2}}U_{2} & = x^{2}U_{1}+\frac{1}{x^{3}}U_{3}
+\lambda_{2}U_{1}-\lambda_{1}U_{2}, &\qquad {\sc D}_{U_{2}}U_{3}
& = -\frac{1}{x^{3}}U_{2} - \frac{1}{2}(x^{3})^{2}U_{1},\\[0.5em]
{\sc D}_{U_{3}}U_{1} & = \frac{1}{2}(x^{3})^{2}U_{2},&\qquad
{\sc D}_{U_{3}}U_{2} & = -\frac{1}{2}(x^{3})^{2}U_{1},\\[0.5em]
{\sc D}_{U_{3}}U_{3} & = 0,
\end{aligned}
\end{equation*}
where $\lambda_{1}$ and $\lambda_{2}$ are real numbers. Since ${\sc B}^{\sc d}={\sc D}-\nabla$, then we obtain
\begin{eqnarray*}
{\sc B}^{\sc d}_{U_{1}}U_{1} &=&\lambda_{1}U_{1}+\lambda_{2}U_{2},\qquad {\sc B}^{\sc d}_{U_{1}}U_{2} = -\lambda_{1}U_{1}-\lambda_{2}U_{2},\\
{\sc B}^{\sc d}_{U_{2}}U_{2} &=&\lambda_{2}U_{1}-\lambda_{1}U_{2},\qquad
{\sc B}^{\sc d}_{U_{1}}\xi = {\sc B}^{\sc d}_{U_{2}}\xi = {\sc B}^{\sc d}_{U_{3}}\xi =0.
\end{eqnarray*}
Then by Theorem \ref{th-tssf-sec-cur-inequality}, we see that $({\sc D},\varphi,\xi,\eta,g)$ is a trans-Sasakian statistical structure of pointwise constant statistical $\varphi$-sectional curvature $-\left(\frac{3}{4}(x^{3})^{4}+\frac{1}{(x^{3})^{2}}+(x^{2})^{2}\right) -2\lambda_{1}^{2}-2\lambda_{2}^{2}$.
\end{ex-new}

\begin{ex-new}
We reconsider Example~\ref{ex-tsm-expo}, in which $M$ is a trans-Sasakian space form with pointwise constant $\varphi$-sectional curvature $-\frac{3}{4}e^{4x^{3}}-1$. Define an affine connection ${\sc D}$ on $M$ as
\begin{equation*}
\begin{aligned}
{\sc D}_{U_{1}}U_{1} & = U_{3}+\frac{a}{2}(U_{1}+U_{2}), &\qquad {\sc D}_{U_{1}}U_{2} & = -\frac{1}{2}e^{2x^{3}}U_{3} + \frac{a}{2}(U_{1}-U_{2}),\\[0.5em]
{\sc D}_{U_{1}}U_{3} & = -U_{1}+\frac{1}{2}e^{2x^{3}}U_{2}, &\qquad {\sc D}_{U_{2}}U_{1} &= \frac{1}{2}e^{2x^{3}}U_{3} - x^{2}e^{x^{3}}U_{2} + \frac{a}{2}(U_{1}-U_{2}),\\[0.5em]
{\sc D}_{U_{2}}U_{2} & = U_{3} + x^{2}e^{x^{3}}U_{1} - \frac{a}{2}(U_{1}+U_{2}),&\qquad{\sc D}_{U_{2}}U_{3} &= -U_{2} - \frac{1}{2}e^{2x^{3}}U_{1}, \\[0.5em]
{\sc D}_{U_{3}}U_{1} & = \frac{1}{2}e^{2x^{3}}U_{2},&\qquad {\sc D}_{U_{3}}U_{2} & = -\frac{1}{2}e^{2x^{3}}U_{1}, &\\[0.5em]
{\sc D}_{U_{3}}U_{3} &=0,
\end{aligned}
\end{equation*}
where $a$ is a real number. Since ${\sc B}^{\sc d}={\sc D} -\nabla$, then we obtain
\begin{eqnarray*}
{\sc B}^{\sc d}_{U_{1}}U_{1} &=&\frac{a }{2}\left(U_{1}+U_{2}\right) ,
\qquad {\sc B}^{\sc d}_{U_{1}}U_{2} =\frac{a}{2}\left(U_{1}-U_{2}\right),\\
{\sc B}^{\sc d}_{U_{2}}U_{2} &=&-\frac{a }{2}\left(U_{1}+U_{2}\right) ,
\qquad {\sc B}^{\sc d}_{U_{1}}\xi = {\sc B}^{\sc d}_{U_{2}}\xi = {\sc B}^{\sc d}_{U_{3}}\xi =0.
\end{eqnarray*}
Then by Theorem \ref{th-tssf-sec-cur-inequality}, we see that $({\sc D},\varphi,\xi,\eta,g)$ is a trans-Sasakian statistical structure of
pointwise constant statistical $\varphi$-sectional curvature $-\frac{3}{4}e^{4x^{3}}-1-a ^{2}$.
\end{ex-new}

\vspace{3mm}
\noindent {\bf Acknowledgement.} The first author gratefully acknowledges the Department of Science and Technology (DST), Government of India, for 
financial support through the INSPIRE Fellowship DST/INSPIRE/03/2022/
002296. The second author is thankful to the Banaras Hindu University for the ``Incentive Grant" under the IoE Scheme (BHU), Ref.\ No.\ R/Dev/D/IoE/Incentive (Phase-IV)/2024-25/82489.

\end{document}